\documentclass[12pt,twoside]{article}

\DeclareSymbolFont{AMSb}{U}{msb}{m}{n}
\DeclareSymbolFontAlphabet{\Bbb}{AMSb}
\usepackage{epsfig}
\newcommand{\QED}{{\hspace*{\fill}\rule{2mm}{2mm}\linebreak}}

\newcommand{\be}{\begin{equation}}
\newcommand{\ee}{\end{equation}}
\newcommand{\ba}[1]{\begin{array}{#1}}
\newcommand{\ea}{\end{array}}
\newcommand{\demo}[1]{\vspace{2mm}\par\noindent{\it #1.\/}}

\newcommand{\R}{{\Bbb{R}}}
\newcommand{\RR}{{{\Bbb{R}}^{3}}}

\newcommand{\Z}{{\Bbb{Z}}}
\newcommand{\ZZ}{{{\Bbb{Z}}^{3}}}
\newcommand{\Zd}{{{\Bbb{Z}}^{d}}}
\newcommand{\Ld}{{{\Bbb{L}}}}
\newcommand{\Ed}{{{\Bbb{E}}}}
\newcommand{\Hd}{{{\Bbb{H}}}}
\newcommand{\Bigmid}{\,\Big|\,}

\newcommand{\Prob}{\Bbb{P}}

\newcommand{\GG}{\mathbf{G}}
\newcommand{\hsep}{H_{\mbox{\rm\scriptsize sep}}}

\newcommand{\la}{\left\langle}
\newcommand{\ra}{\right\rangle}

\newcommand{\rc}{{\rm c}}
\newcommand{\ol}[1]{\overline{#1}}
\newcommand{\ul}[1]{\underline{#1}}
\newcommand{\rd}{{\rm d}}
\newcommand{\es}{\emptyset}
\newcommand{\Dv}{\Delta_{\rm v}}
\newcommand{\De}{\Delta_{\rm e}}
\newcommand{\ve}{{\mbox{\rm\scriptsize v}}}

\newcommand{\smbox}[1]{\mbox{  #1{} }}

\newcommand{\La}{\Lambda}
\newcommand{\Om}{\Omega}

\newcommand{\marginal}[1]{\strut\setbox0=%
      \vtop{\hsize=15mm
            \parindent=0pt\baselineskip=8pt
             \raggedright\rightskip=1em plus 2em\hfuzz=.5em\tolerance=9000
             \overfullrule=0pt
             #1}%
      \dimen0=\ht0
      \advance\dimen0 by \dp\strutbox
      \ht0=0pt\dp0=\dimen0
      \vadjust{\kern-\dimen0\moveleft\wd0\box0}%
      \ignorespaces}       

\newtheorem{theo}{Theorem}
\newtheorem{pro}{Proposition}
\newtheorem{lem}[pro]{Lemma}

\newcounter{mycount}
\newenvironment{mylist}{\begin{list}{\rm(\roman{mycount})}%
   {\usecounter{mycount}\labelwidth=1cm\itemsep 0pt}}{\end{list}}

\newcommand{\pc}{p_{\rm c}}
\newcommand{\bond}[1]{\langle#1\rangle}
\newcommand{\inter}{\mbox{\rm int}}
\newcommand{\ins}{\mbox{\rm ins}}
\newcommand{\out}{\mbox{\rm out}}
\newcommand{\odelta}{{\overline \delta}}
\newcommand{\sdelta}{{\delta^\ast}}

\def\0sim{\stackrel{0}{\sim}}
\def\1sim{\stackrel{1}{\sim}}
\newcommand{\DLM}{{\cal D}_{L, M}}
\newcommand{\ILM}{{\cal I}_{L, M}}
\newcommand{\OLM}{\Omega_{L,M}^\mu}
\newcommand{\DL}{{\cal D}_L}
\newcommand{\IL}{{\cal I}_L}
\newcommand{\tfrac}{\textstyle\frac}
\newcommand{\lest}{\le_{\rm st}}
\newcommand{\gest}{\ge_{\rm st}}
\newcommand{\q}{\quad}
\newcommand{\lra}{\leftrightarrow}
\newcommand{\de}{\delta}
\newcommand{\of}{f^o}

\begin{document}

\begin{titlepage}

\vspace{5mm}
\hfill 13 September 2001 

\begin{center}
 {\large {\bf  Rigidity of the interface in percolation\\ and random-cluster models}}
\\[10mm]

Guy Gielis\footnote[1]{King's College Research Centre,
Cambridge CB2 1ST,
and
Statistical Laboratory,
Centre for Mathematical Sciences,
Cambridge CB3 0WB, UK. E-mail:
guy.gielis@barco.com}
\\[2mm]

Geoffrey Grimmett\footnote[2]{Statistical Laboratory, 
Centre for Mathematical Sciences, 
Cambridge CB3 0WB, UK. E-mail:
g.r.grimmett@statslab.cam.ac.uk, http://www.statslab.cam.ac.uk/$\sim$grg/}
\\[10mm]

\end{center}

\noindent
{\bf Abstract:}
We study conditioned random-cluster measures with edge-parameter $p$ and cluster-weighting factor
$q$ satisfying $q\ge 1$. The conditioning corresponds to mixed boundary conditions for a spin
model. Interfaces may be defined in the sense of Dobrushin, and these are proved
to be `rigid' in the thermodynamic limit, in three dimensions and for sufficiently large values
of $p$. This implies the existence of non-translation-invariant 
(conditioned) random-cluster 
measures in three dimensions. The results are valid 
in the special case $q=1$, thus indicating
a property of three-dimensional percolation not previously noted.

\vspace{3mm}
\noindent
{\bf Keywords:} Random-cluster model, percolation, Ising model,
Potts model, interface, Dobrushin boundary condition.

\vspace{3mm}
\noindent
{\bf Mathematics Subject Classification (2000):} 60K35, 82B20, 82B43.

\end{titlepage}


\section*{1. Introduction}
Dobrushin's proof \cite{Do} of the existence of non-translation-invariant Gibbs states for the
three-dimensional Ising model was the starting point for the study of interfaces
in disordered spin systems. We show in the current paper that such results are valid for all
ferromagnetic random-cluster models on $\ZZ$, 
including percolation. This generalization
of Dobrushin's theorem is achieved by defining a family of conditioned measures, and by showing
the stiffness of the ensuing interface.

The random-cluster model has since its introduction 
\cite{Fo1, Fo2, FK} around 1970 provided a beautiful mechanism
for the study of Ising and Potts models, as well as being an object worthy of
study in its own right. Many (but not all) central
results for ferromagnetic Ising/Potts systems are best proved in the context of
random-cluster models; the stochastic geometry of such models may be exploited
the better to understand the behaviour of correlations in the original system.
The spectrum of random-cluster models extends to percolation (and beyond), and one
sees thus that percolative techniques have direct application to Ising and Potts models.
The reader is referred to \cite {Gr0} for more information concerning the history of random-cluster
models, and to \cite {ACCN, BGK, CP,  GHM, Gr1} for examples of them in action.

The question addressed here concerns the stiffness of interfaces. In the case of the Ising model,
Dobrushin introduced the boundary condition on the box $\Lambda=[-L,L]^3$ having $+1$ on the upper
half of the boundary and $-1$ on its complement. He then studied the interface separating the two regions
behaving respectively as the $+1$ phase and the $-1$ phase. He showed 
for sufficiently low temperatures that this interface deviates only locally from 
the horizontal plane through
the equator of the box. This effect is seen in all dimensions of three or more, but not
in two dimensions, for which case the interface may be thought of as
a line with
Gaussian fluctuations (see \cite {GH, Hr}).

This problem may be cast in the more general setting of the random-cluster model on the box
$\Lambda$ subject to the following boundary condition and to a certain conditioning.
The vertices on the upper hemisphere of $\Lambda$ are wired together into a single composite
vertex labelled $\Lambda^+$. The vertices on the complement of the upper
hemisphere are wired into a single composite vertex labelled $\Lambda^-$. Let $\cal D$
be the event that no open path of $\Lambda$ exists 
joining $\Lambda^-$ to $\Lambda^+$, and let $\phi_\La$ be the
random-cluster measure on $\Lambda$ with edge-parameter $p$ and cluster-weighting
factor $q$, with the above boundary condition {\it and conditioned on the event\/} $\cal D$.
It is a geometrical fact that there
exists an interface separating two regions of $\Lambda$, each of which is in the wired phase.
It follows by the results of \cite {Do} that, when $q=2$ and $p$ is sufficiently
large,  this interface deviates only locally from the horizontal plane through the equator
of $\Lambda$. The purpose of this paper is to prove that this is so for all $q\ge 1$
and all sufficiently large $p$. In doing so we shall work directly with the
random-cluster model. The geometry of the interfaces for this model is notably
different from that of a spin model since the configurations
are indexed by edges rather than by vertices, and this leads to some new difficulties.

Extensions of our results to dimensions $d$ satisfying $d\ge 4$ are, to quote from \cite{Do},
``obvious'', though the proofs may involve
some extra complications. It is striking that our results are valid for 
high-density percolation on $\ZZ$, being the case $q=1$.
That is, conditional on the existence of a surface (suitably defined) of dual plaquettes
spanning the equator of $\Lambda$, this surface deviates only locally from the flat plane.
A corresponding question for supercritical percolation in two dimensions has been studied in
depth in \cite{CCC}, where it is shown effectively that the (one-dimensional) interface
converges when re-scaled to a Brownian bridge.

We have spoken above of interfaces which `deviate only locally' from a plane, and we shall make this
expression more rigorous in Section 9, where our principal
Theorem 2 is presented. We include in Section 3
a weaker version of Theorem 2 which does not make use of the notation
developed later in the work.

Our theorems are proved under the assumption that $q\ge 1$ and $p$ is 
sufficiently large. It is a major open
question to determine whether or not such results are valid under the weaker assumption that
$p$ exceeds the critical value $\pc(q)$ of the random-cluster model with cluster-weighting
factor $q$ (see \cite{Gr2}).  The answer may be expected to depend
on the value of $q$ and the number $d$ of dimensions. Since the percolation measure
(when $q=1$) is a conditioned product measure, it may be possible as
with other problems to gain insight into the existence or not of a `roughening transition'
by concentrating on the special case of percolation.
It is of interest that much of the argument 
of this paper is valid also when $q<1$ and $p$ is
sufficiently large, but we shall not specify the details.
Also, it may be possible to extend some of the conclusions of this paper
to measures with certain other boundary conditions, such
as that generated with free boundary conditions
and conditioned on $\cal D$, but
we shall not pursue this here.

As described above, the measures studied here are obtained by
conditioning on a certain event $\cal D$. When $p$ is large, $\cal D$ has probability
of order $\exp(-\alpha L^2)$ where $\alpha=\alpha(p,q)$, and thus we are in the realm of the 
large-deviation theory of the process. See \cite{CP, DP}.

We introduce random-cluster measures in the next section, 
followed by a summary of our main results in Section 3.
Necessary properties of random-cluster measures are developed in Section 4.
Interfaces are defined in Section 5, where we prove some 
geometrical lemmas of independent interest which we believe
will find applications elsewhere. 
In Section~6 we study the probability of having a configuration
that is compatible with a given interface, under the appropriate
conditioned measure. We present in
Section~7 a microscopic geometrical description of the
random-cluster interfaces using a terminology
based on that introduced for the Ising model in \cite{Do}. 
This is followed in Section 8 by an exponential bound
for the probability of finding local perturbations
of a flat interface, and in Section 9 by the statement and proof
of our main theorems.

The methods of this paper are inspired by those of \cite {Do} subject to
some serious variations. Dobrushin \cite{Do} studied the Ising model,
and his arguments were later simplified in part by van Beijeren \cite{Be}.
We have been unable to extend the methods of \cite {Be}, which may
be special to the Ising model. Related results may be found in
\cite{BLOP, BLP, DMN, HKZ, HZ} and the references therein. We have found the first
of the latter references to be particularly useful in the present work.
It should be noted
that, in order to study interfaces for spin systems rigorously, certain lemmas concerning their
geometry are required; see \cite{HKZ, KLMR} for example.

The Pirogov--Sinai theory of contours has enabled 
(\cite{KLMR, LMMRS, MMRS}) a study of Potts models
and random-cluster models for large $q$, when $p=\pc(q)$, the critical point.
It seems now to be accepted that the random-cluster model is especially
well adapted to the study of contours and interfaces. However, it appears
that certain pivotal facts, implicit in earlier work, and concerning
the relationship between interfaces and random-cluster measures, have never been proven.
Specifically, certain key results in three dimensions concerning the `external
boundary' of a set of connected edges, and the `internal boundary' of
a cavity of plaquettes of $\ZZ$, are missing from the literature.
These are akin to the well known fact, proved in \cite {Ke82},
that the external boundary of a finite cluster of $\Z^2$ contains, in its dual representation, 
a circuit separating the cluster from infinity. One
of the targets of the current paper is to state and prove the necessary
geometrical facts; see Propositions \ref{P} and \ref{Q}.

Since finishing this work, we have received the preprint \cite{CK},
which uses Pirogov--Sinai theory to study
the rigidity of interfaces for sufficiently large $q$ and with $p$ equal
to the critical point of the random-cluster model. It is proved there that
there is a rigid interface at a first-order transition for large $q$,
with the boundary condition a mixture of the wired
and the free.


\section*{2. Conditioned random-cluster measures}
Let $\ZZ$ be the set of all vectors $x=(x_1,x_2,x_3)$ of integers, termed {\it vertices\/}, and let
$$
| x-y| = \sum_{i=1}^3 |x_i-y_i|,\ \| x-y\| = \max\Big\{|x_i-y_i|: 1\le i\le 3\Big\}
\quad\mbox{for } x,y\in\RR.
$$
We place an {\it edge\/} between vertices $x$ and $y$ if and only if $|x-y|=1$, and we
denote by $\Ld=(\ZZ,\Ed)$ the resulting lattice. 
We write $x \sim y$ if $|x-y|=1$, and we write $\langle x,y\rangle$
for the corresponding edge. We sometimes think of the edge $e=\langle x,y\rangle$
as the closed straight-line segment with endpoints $x$ and $y$.
For $E\subseteq \Ed$, we write $V(E)$ for the set of vertices in $\ZZ$ that belong to at 
least one of the edges in $E$. We shall sometimes abuse notation by referring to the graph
$(V(E),E)$ as the graph $E$. The $L^\infty$ distance between two edges $e_1$, $e_2$
is defined to be
the distance between their centres, and is denoted $\|e_1,e_2\|$.

A {\it path\/} in a subgraph $G=(V,E)$ 
of $\Ld$ is an alternating set of distinct vertices and
bonds $x=z_0, \bond{z_0,z_1}, z_1, \ldots, \bond{z_{n-1},z_n}, z_n=y$ using only edges
$\bond{z_i,z_{i+1}}\in E$. Such a path is said to connect $x$ and $y$
and to have length $n$.
The graph $G$ is called {\it connected\/} if every pair of vertices is connected
by some  path. A {\it connected component\/} of $G$ is a 
maximal connected subgraph of $G$.
We shall occasionally speak of a set $\La\subseteq \ZZ$ of vertices as being {\it connected\/}, 
by which we mean that $\La$ induces a  connected subgraph of $\Ld$. 

For $x\in\ZZ$, we denote by $\tau_x:\ZZ\to\ZZ$
the translate given by $\tau_x(y)=x+y$. The translate $\tau_x$
acts on edges and subgraphs of $\Ld$ in the natural way.
For sets $A$, $B$ of edges or vertices of $\Ld$, we write
$A \simeq B$ if $B=\tau_x A$ for some $x\in\ZZ$. Note that
two edges $e$, $f$ satisfy $\{e\}\simeq \{f\}$ if and only if
they are parallel (in which case we write $e\simeq f$).

We write $S^\rc$ for the complement of a set $S$.
The {\it upper\/} and {\it lower boundaries\/} of a set $\La$ of vertices
are defined as
$$
\ba{rl}
 \partial^+\La&=\big\{x\in \La^\rc: x_3>0, \ x\sim z\mbox{ for some } z\in \La  \big\}, \\
 \partial^-\La&=\big\{x\in \La^\rc: x_3\le 0, \ x\sim z\mbox{ for some } z\in \La  \big\}, 
\ea
$$
and the {\it boundary\/} of $\La$ is denoted
$\partial \La=\partial^+\La\cup\partial^-\La$. 
For positive integers $L$, $M$ we define the box $\La_{L,M}=[-L,L]^2\times[-M,M]$,
and write $E_{L,M}$ for the set of
all edges having at least one endvertex in $\La_{L,M}$. 
[We abuse notation here and later, and should write 
$\La_{L,M}=([-L,L]^2\times[-M,M])\cap \ZZ$.]
We write $Q_L=\La_{L,L}$, the cube of side-length $2L$,
and $\La_L=[-L,L]^2\times \Z$, an infinite cylinder.

The configuration space of the random-cluster model on $\Ld$
is the set $\Omega=\{0,1\}^\Ed$, which we endow with the $\sigma$-field
$\cal F$ generated by the finite-dimensional cylinders. 
A configuration $\omega\in\Omega$ assigns to each edge $e$  the 
value $0$  or $1$; we call the edge $e$ {\it open\/} (in $\omega$) if
$\omega(e)=1$, and {\it closed\/} otherwise. A set of edges (for example, a path) is
called {\it open\/} if all the edges therein are open.
For $\omega\in\Omega$,
we write $x\lra y$ if there exists an open path connecting the vertices $x$ and $y$,
and $x\lra A$ if there exists $y\in A$ such that $x \lra y$.
Each $\omega\in\Omega$ is in one--one correspondence with its set  
$\eta(\omega)=\{e\in \Ed:\omega(e)=1\}$ of open edges. We write $\eta_x(\omega)$
for the set of edges in the connected component of the graph
$(\ZZ,\eta(\omega))$ containing the vertex $x$.
The configuration
which assigns 1 (respectively 0) to every edge is denoted 1 (respectively 0).

Let $E$ be a finite subset of $\Ed$ and let $V=V(E)$, and suppose that $0\le p\le 1$ and $q>0$.
The usual way (see \cite{Gr2}) of defining a random-cluster measure with parameters $p$, $q$
on the graph $G=(V,E)$ with
boundary condition $\zeta$ ($\in\Omega$) is via the formula
$$
   \phi_{G,p,q}^{\zeta}(\omega)=\frac{1}{Z^\zeta_{G,p,q}}\left\{\prod_{e\in E}p^{\omega(e)}(1-p)^{1-\omega(e)}\right\}
        q^{k_G(\omega)}
     I\{\omega(f)=\zeta(f) \mbox{\ if\ } f\notin E\},
$$
defined for all $\omega\in\Omega$.
Here, $k_G(\omega)$ is the number of connected components 
in the  graph $(\ZZ, \eta(\omega) )$ having at least one vertex belonging to
$V$, 
\be
\label{part}
Z^\zeta_{G,p,q}=\sum_{\omega\in\Omega}\left\{\prod_{e\in E}p^{\omega(e)}(1-p)^{1-\omega(e)}\right\}
        q^{k_G(\omega)}
     I\{\omega(f)=\zeta(f) \mbox{\ if\ } f\notin E\}
\ee 
is the normalizing partition function, and $I\{H\}$ is the indicator function 
of the event $H$. We shall write $k(\omega)$ for the total number of
connected components of $(\ZZ,\eta(\omega))$.

We shall be particularly concerned with the case $E=E_{L,M}$ and with a boundary condition
$\mu$ corresponding to the mixed `Dobrushin boundary'
of \cite {Do}. To this end, we let $\mu$ be given by
\be
\label{iota}
\mu(e)= \left\{\begin{array}{ll}
  0   & \smbox{if} e=\bond{x,y} \mbox{ for some } x=(x_1,x_2,0) 
                                 \smbox{and} y= (x_1,x_2,1),\\
  1  & \smbox{otherwise.}
\end{array}\right.
\ee
We let $\Omega^\mu_{L,M}$ be the set of all configurations $\omega\in\Omega$
such that $\omega(f)=\mu(f)$ if $f\notin E_{L,M}$. We define ${\cal I}_{L,M}$ 
to be the event that
there exists no open path connecting a vertex of
$\partial^+ \La_{L,M}$ to a vertex of  $\partial^- \La_{L,M}$. 
Let $\ol\phi^\mu_{\La_{L,M},p,q}$ denote the measure $\phi_{G,p,q}^\mu$ conditioned
on the event ${\cal I}_{L,M}$. See Figure 1.

\begin{figure}
\centering
\epsfig{file=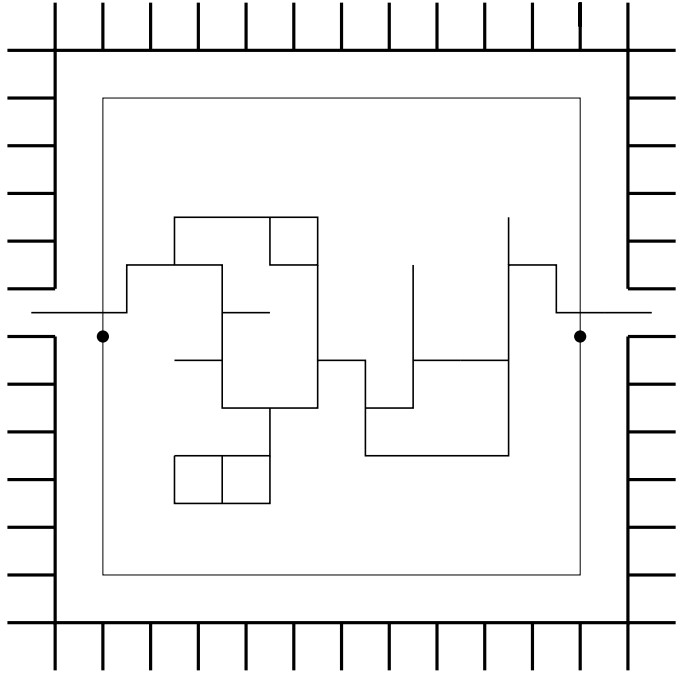,width=9cm,height=9cm}
\caption{The box $\La_{L,M}$. The heavy black edges
are those given by the boundary condition 
$\mu$, and there is a two-dimensional sketch of the interface $\Delta$.}
\label{slide}
\end{figure}

The measure $\ol\phi^\mu_{\La_{L,M},p,q}$ is only one of many such conditioned measures. 
Let $E$ be a finite subset of $\Ed$, let $V=V(E)$, and write $G=(V,E)$ as usual.
In a more
general formulation, we take some boundary condition $\zeta$, and we consider
the set $\cal C(\zeta)$ of open components of $\zeta$ in the graph 
obtained from $\Ld$ by removing
both $E$ and all vertices adjacent to no edge in $E^\rc$. 
Let $S$ be some set of labels, let $l:{\cal C}(\zeta)\to S$, and call
$l(C)$ the {\it label\/} of $C\in {\cal C}$.
We now consider the measure $\phi^\zeta_{G,p,q}$
conditioned on the event that no open path exists joining two vertices lying in components
of $\cal C(\zeta)$ having different labels, and we denote this new measure
by $\ol\phi^{\zeta,l}_{G,p,q}$. The case above arises when $E=E_{L,M}$ and $\zeta=\mu$, 
(note that $|{\cal C}(\mu)|=2$),
and the two members of $\cal C(\mu)$ have different labels. 


\section*{3. Summary of main results}
\noindent
We summarise our main results as follows. The complete
form of our main theorems appear with proofs in Section 9,
using notation developed in the course of the work.

Many of our calculations concern the box $\La_{L,M}$ and the
measure $\ol\phi^\mu_{\La_{L,M},p,q}$. We choose however to express our conclusions in terms
of the infinite cylinder $\La_L=\La_{L,\infty}$ and the weak limit
$\ol\phi_{L,p,q}=\lim_{M\to\infty}\ol\phi^\mu_{\La_{L,M},p,q}$, which
is shown in Lemma \ref{limitphi} to exist.

We show in Proposition \ref{intuition} that, on the event $\ILM\cap \Omega^\mu_{L,M}$, there
exists an `interface which spans the equator' of $\La_{L,M}$.
(By the equator, we mean the circuit of $\La_{L,M}\setminus\La_{L-1,M}$
comprising all vertices $x$ with $x_3=\frac12$.)
Much of this paper is devoted to understanding the geometry of such an interface.
We shall see in Theorem \ref{theorem1} that, in the limit as $M\to\infty$
and for sufficiently large $p$, this interface deviates, $\ol\phi_{L,p,q}$-almost surely,
only locally from the flat plane through the equator of $\La_L$.
Indeed, the spatial density of such deviations approaches zero
as $p$ approaches $1$.
As a concrete application we present the following theorem, which we note
to be a substantial weakening of Theorem \ref{theorem1}
in Section 9.

\begin{theo}
\label{theorem2}
Let $q\geq 1$. For all $\epsilon>0$ there exists 
$\widehat p=\widehat p(\epsilon)<1$ such that, if $p>\widehat p$,
\be
\label{theorem1res}
  \ol\phi_{L,p,q}\Big(x\lra\partial^- \La_L\Big)>1-\epsilon,
\quad \ol\phi_{L,p,q}\Big(x+(0,0,1)\lra\partial^+ \La_L\Big)>1-\epsilon,
\ee
for all $L\ge 1$ and every $x=(x_1,x_2,0)\in \La_{L}$.  
\end{theo}

We have no proof that the sequence $\{\ol\phi_{L,p,q}:L\ge 1\}$ converges
weakly as $L\to\infty$, but, by the usual compactness argument, there must
exist weak limits of the sequence. It is a consequence 
of our main Theorem \ref{theorem1} that, for
sufficiently large $p$, any such weak limit is non-translation-invariant.
By making use of the relationship between
random-cluster models and Potts models (see \cite{ACCN, Gr0} and
the references therein),
one obtains thereby a generalization
of the theorem of Dobrushin \cite{Do} to include
percolation and Potts models.
We return to this point in Section 9, where it is shown in addition that
there exists a geometric bound, uniform in $L$, on the tail
of the displacement of the interface from the flat plane.

It would be interesting to know more of the random field
defined by the locations where the interface coincides with the flat plane through the equator.
Could it be the case that this field dominates (stochastically)
a percolation process of some density $\rho(p)$, where $\rho(p)\to 1$ as $p \to 1$?
Such a proposal is supported, in the special
case when $q=1$ at least, by the new correlation/connection inequality
of \cite{BK}.

\section*{4. Properties of random-cluster measures}
There follow some general lemmas concerning random-cluster measures.
The first of these contains the comparison inequalities
of Fortuin and Kasteleyn. There is  a partial order on $\Omega$
given by $\omega_1\le \omega_2$ if and only if $\omega_1(e)\le \omega_2(e)$
for all $e\in\Ed$.  A function $h:\Omega\to\R$ is called {\it increasing\/}
if it is increasing with respect to this partial order.
Given two probability measures $P_1$, $P_2$ on $(\Omega,{\cal F})$, we write $P_1\lest P_2$ if
$\int h \,dP_1 \le \int h\,dP_2$
for all bounded measurable increasing functions $h$.

\begin{lem}
\label{comp}
Let $E$ be a finite subset of $\Ed$, and $G=(V,E)$ where $V=V(E)$.
For any $\zeta\in\Omega$,
we have that
\be
\label{compeq}
\ba{rl}
\phi^\zeta_{G,p',q'}\lest \phi^\zeta_{G,p,q}\q&\mbox{if}\q p'\leq p,\ q'\geq q,\ q'\geq1,\\
\phi^\zeta_{G,p',q'}\gest\phi^\zeta_{G,p,q}\q&\mbox{if}\q\displaystyle\frac{p'}{q'(1-p')}
\geq\frac{p}{ q(1-p)},\ q'\geq q,\ q'\geq 1.
\ea
\ee
\end{lem}

See \cite{ACCN, Gr2} for a proof of these standard inequalities. Our second lemma is
a formula for the partition function in terms of the edge densities. 
For $e\in \Ed$, we write $J_e$ for the event that $e$ is open.

\begin{lem}
\label{logint}
Let $E$ be a finite subset of $\Ed$, and $G=(V,E)$ where $V=V(E)$.
For any $\zeta\in\Omega$,
we have that
\be
\log Z^\zeta_{G,p,q}= k_G(\zeta^1)\log q +\sum_{e\in E} g^\zeta_{G,p,q}(e),
\ee
where $\zeta^1$ is the configuration obtained from $\zeta$
by making every edge in $E$ open,  and
\be
\label{logz}
g^\zeta_{G,p,q}(e)=\int_p^1 
\left[\displaystyle \frac{r-\phi^\zeta_{G,r,q}(J_e)}{r(1-r)}\right]\,dr.
\ee
\end{lem}

\demo{Proof}
We differentiate $\log Z^\zeta_{G,r,q}$ with respect to $r$, 
as in \cite {Gr2}, p.\ 1479, to obtain
that
$$
\displaystyle \frac{d}{dr}\log Z^\zeta_{G,r,q}=
\sum_{e\in E}\frac{\phi^\zeta_{G,r,q}(J_e)-r}{r(1-r)}.
$$
This we integrate from $p$ to 1, and note that $\log Z^\zeta_{G,1,q}=k_G(\zeta^1)\log q$.
\QED

Let $q \ge 1$.
We have by Lemma 1 that $\phi^\zeta_{G,r',1}\lest \phi^\zeta_{G,r,q}\lest 
\phi^\zeta_{G,r,1}$ where $r'=r/(r+(1-r)q)$, and hence
$$
\displaystyle\frac r{r+(1-r)q} \le \phi^\zeta_{G,r,q}(J_e)\le r.
$$
By substitution into (\ref{logz}), 
\be
\label{bndg}
0 \le g^\zeta_{G,p,q}(e) \le \int_p^1 (q-1)\,dr = (1-p)(q-1)\q\mbox{for } e\in E,
\ee
uniformly in the choice of $E$ and $\zeta$.
The above inequalities are reversed if $q< 1$.

We recall for the next lemma that $Q_n=\La_{n,n}$, and, for $e\in\Ed$,
we write $Q_n(e)=e+Q_n$, the set of translates of 
the endvertices of $e$ by 
vectors in $\La_{n,n}$.

\begin{lem}
\label{expbnd}
Let $q\ge 1$. There exists $p^\ast=p^\ast(q)<1$ and a
constant $\alpha>0$
such that the following holds. Let $E_1$ and $E_2$ be finite edge sets
such that $e\in E_1\cap E_2$, and let $n\ge 1$ be 
such that $E_1\cap Q_n(e)=E_2\cap Q_n(e)$. If $p>p^\ast$,
$$
|g^1_{G_1,p,q}(e)-g^1_{G_2,p,q}(e)| \le e^{-\alpha n},
$$
where $G_i=(V(E_i),E_i)$.
\end{lem}

\demo{Proof}
Let $L_e$ be the event that the endvertices of the edge $e$ are joined by 
an open path which does not use $e$ itself. It is an elementary argument, using
equation (3.10) of \cite{Gr2}, that
$$
\displaystyle \frac{r-\phi^1_{G,r,q}(J_e)}{r(1-r)} = \frac{(q-1)(1-\phi^1_{G,r,q}(L_e))}
{r+(1-r)q},
$$
whence
\be
\label{g1g2}
|g^1_{G_1,p,q}(e)-g^1_{G_2,p,q}(e)| \le \int_p^1 \frac{(q-1)}{r+(1-r)q}
|\phi^1_{G_1,r,q}(L_e) -\phi^1_{G_2,r,q}(L_e) |\,dr.
\ee

Fix $n\ge 1$.
We shall now follow an argument of \cite{Gr2}, pp.\ 1486--1487, and
\cite{Ke}, pp.\ 138--152, of which we give some details next.
Let $\cal L$ be derived from $\Ld$ by adding edges between any
pair $x$, $y$ of vertices with $\|x-y\|=1$. For $\omega\in\Omega$, we
call a vertex $x$ {\it white\/} if $\omega(e)=1$ for all $e$ incident with
$x$ in $\Ld$, and {\it black\/} otherwise.
Let $V$ be the set of vertices which are incident in $\Ld$ to edges of
both $Q_n(e)$ and its complement. We define $B$ as the union of $V$
together with all vertices $x_0\in\ZZ$ for which there exists a path
$x_0,x_1,\dots,x_m$ of $\cal L$ such that $x_0,x_1,\dots,x_{m-1}\notin V$,
$x_m\in V$, and $x_0,x_1,\dots,x_{m-1}$ are black. Let $K_n$ be the
event that there exists no $x\in B$ such that $\| x-z\| \le 10$, say,
where $z$ is the centre of $e$.
Using (5.17)--(5.18) of
\cite{Gr2}, together with estimates at the beginning of the proof
of Lemma (2.24) of \cite{Ke}, we find that  
\be
\label{GK}
\phi^0_{Q_n(e),r,q}(K_n)\ge 1-c^n(1-\pi)^{en}
\ee
where $c$ and $e$ are absolute positive constants, and
$\pi=r/(r+(1-r)q)$. Since $K_n$ is an increasing event,
we deduce that
\be
\label{GK2}
\phi^1_{G_1,r,q}(K_n)\ge 1-c^n(1-\pi)^{en}.
\ee
Let $H=E_1\cap Q_n(e)$. It follows by the arguments of \cite{Gr2}, p.\ 1487, and
by coupling, that
$$
0\le \phi^1_{H,r,q}(L_e)-\phi^1_{G_1,r,q}(L_e) \le 1-\phi^1_{G_1,r,q}(K_n).
$$
The claim then follows by (\ref{g1g2}), 
(\ref{GK2}),  and the triangle inequality.
\QED

\section*{5. Interfaces and geometry}
We shall have much recourse to the dual of the random-cluster model, being
a probability measure on the set of `plaquettes' of the dual lattice $\Ld_\rd$
obtained by shifting the vertices and edges of $\Ld$ through the vector
$(\frac12,\frac12,\frac12)$
(see \cite{ACCFR, GrH}).
A {\it plaquette\/} of $\Ld_\rd$ is a (topologically) closed
unit square of $\RR$ with corners lying in $\ZZ+(\frac12,\frac12,\frac12)$. 
We denote by $\Hd$ the set of all plaquettes
of $\Ld_\rd$.
The straight line segment joining the
vertices of an edge $\bond{x,y}$ passes 
through the middle of exactly one plaquette, denoted 
$h(\bond{x,y})$, which we call the {\it dual\/} 
plaquette of $\bond{x,y}$. We 
declare this plaquette {\it open\/} (respectively 
{\it closed\/}) if $\bond{x,y}$ is closed (respectively open). 
The plaquette $h(\bond{x,y})$ is called {\it horizontal\/}
if $y=x+(0,0,\pm 1)$, and {\it vertical\/} otherwise.

Two distinct
plaquettes $h_1$ and $h_2$ are called {\it $0$-connected\/}, written $h_1\0sim h_2$ if $h_1\cap
h_2\neq  \es$. They are said to be {\it $1$-connected}, written $h_1\1sim h_2$, if  $h_1\cap h_2$ is
homeomorphic to the unit interval $[0,1]$. A set of plaquettes is called $0$-connected (respectively
$1$-connected) if they are connected when viewed as the vertex-set of a graph with adjacency relation 
$\0sim$ (respectively $\1sim$). The $L^\infty$
distance between two plaquettes $h_1$, $h_2$ is defined
to be the distance between their centres, and is denoted $\|h_1,h_2\|$. 
For any set $H$ of plaquettes,
we write $E(H)$ for the set of edges of $\Ld$ to which they are dual.

We define the  
{\it regular interface\/} as the set $\de_0$ given by
$$
   \delta_0=\Big\{h\in \Hd: h=h(\bond{x,y}) 
           \smbox{for some} x=(x_1, x_2,0) \smbox{and} y=(x_1, x_2,1)\Big\}.
$$
The {\it interface\/} $\Delta(\omega)$ of a configuration 
$\omega \in {\cal I}_{L,M}\cap \Om^\mu_{L,M}$
is defined to be the maximal
$1$-connected set of open plaquettes containing the plaquettes of 
$\delta_0\setminus \{h(e): e\in E_{L,M} \}$. The set of all interfaces is 
\be
\label{D}
   \DLM =\{\Delta(\omega):\omega \in {\cal I}_{L,M}\cap \Om^\mu_{L,M}\}.
\ee    
While it is tempting to think of an interface
as part of a deformed plane, it may in fact
have a much more complex geometry involving cavities
and attachments.
The following proposition, which will be proved later in this section,
confirms that the interfaces in $\DLM$ separate the top of 
$\La_{L,M}$  from  its bottom.


\begin{pro}\label{intuition}
The event ${\cal I}_{L,M}\cap \Om^\mu_{L,M}$ 
is the set of all configurations $\omega\in\Omega^\mu_{L,M}$
for which there exists $\delta\in\DLM$ such that $\omega(e)=0$ whenever $h(e)\in\delta$.
\end{pro}

For $\delta \in \DLM$, we define its {\it extended interface\/} $\odelta$ to be  the set
\be
\label{deltabar}
   \odelta=\delta\cup\{h\in \Hd: h \mbox{ is $1$-connected to some member of } \delta \}. 
\ee
It will be useful to introduce the  `maximal' ($\ol\omega_\delta$) and `minimal'  
($\ul\omega_\delta$) configurations in $\Omega_{L,M}^\mu$ which are compatible with $\delta$:
\be
\label{minmaxconf}
   \ol\omega_\delta(e)=\left\{ \begin{array}{ll}
          0 & \smbox{if} e\in \delta,\\
          1 & \smbox{otherwise,}\end{array}\right.\quad
    \ul\omega_\delta(e)=\left\{\begin{array}{ll}
               \mu(e) & \smbox{if} e\notin E_{L,M},\\
    		     1  & \smbox{if} e\in E_{L,M}\cap (\odelta\setminus \delta),\\
		     0 &  \smbox{otherwise.}\end{array}\right.
\ee

In Section 6, we shall consider interfaces spanning the
equator of the infinite cylinder $\La_L$.

We consider next some geometrical matters. The words `connected' 
and `component' should be interpreted for the moment
in the topological sense. Let $T\subseteq \RR$, and
write $\ol T$ for the closure of $T$ in $\RR$. 
We define the {\it inside\/} $\ins(T)$ of $T$
to be the union of all the bounded connected
components of $\RR\setminus T$; the {\it outside\/} $\out(T)$ is the union of all the unbounded
connected components of   $\RR\setminus T$. The set $T$ is said to {\it separate\/} $\RR$ if
$\RR\setminus T$ has more than one connected component.  
For a set $H\subseteq\Hd$ of plaquettes,  we define the set $[H]\subseteq \RR$ by 
$[H]=\{x\in \RR: x\in h\mbox{ for some  } h\in H\}$.
We call a finite set $H$   of plaquettes a 
{\it splitting set\/} if $[H]$ is $1$-connected in $\RR$ 
and $\RR\setminus [H]$ contains at least one bounded  connected component.

The following two propositions are in a sense dual to one another, and we believe
they will find applications elsewhere. 
The first is an analogue in three dimensions
of Proposition 2.1 of the Appendix of \cite{Ke82}, where two-dimensional
mosaics are considered.


\begin{pro}\label{P}
Let $G=(V,E)$ be a finite connected subgraph of $\Ld$.
There exists a splitting set $Q$ of plaquettes such that\/{\rm :}
\begin{mylist}
\item $V\subseteq \ins([Q])$,
\item every plaquette in $Q$ is dual to 
some edge of $\Ed$ having exactly one endvertex in 
      $V$, 
\item if $W$ is a connected set of vertices such that $V\cap W=\es$,
and there exists an infinite  path on $\Ld$ starting in $W$  which
uses no vertices in  $V$, then $W\subseteq\out([Q])$.      
\end{mylist}
\end{pro}

Let $\de=\{h(e): e\in D\}$ be a 1-connected set of plaquettes, 
and let $\ol\de$ be given as in
(\ref{deltabar}). Consider the graph $(\ZZ,\Ed\setminus D)$, and let $C$
be a connected component of this graph. Let $\Dv C$ be the set of all vertices $v$ in $C$ for
which there exists $w\in\ZZ$ with $h(\bond{v,w})\in\ol\de$, and let
$\De C$ be the set of edges $f$ of $C$ for which $h(f)\in\ol\de\setminus \de$.
Note that edges in $\De C$ have both endvertices belonging to $\Dv C$.

\begin{pro}\label{Q}
For any finite connected component $C$ of the graph $(\Zd,\Ed\setminus D)$,
the graph $(\Dv C, \De C)$ is connected.
\end{pro}

We shall apply this proposition in the following way.
Let $\de\in\DLM$.
Consider the connected components of the graph
$(\ZZ,\eta(\ol\omega_\delta))$, and denote
these components as $(S_\delta^i, U_\delta^i)$, $i=1,2,\ldots, K_\delta$, where
$K_\delta=k(\ol\omega_\delta)$. Note that
$U_\delta^i$ is empty whenever $S_\delta^i$ is a singleton. 
We define $W(\de)$ as the set of edges in
$E_{L,M}\setminus\{e\in \Ed: h(e)\in \odelta \}$.

Let $\omega\in\ILM\cap \Omega^\mu_{L,M}$ 
be such that $\Delta(\omega)=\de$. It must
be the case that
\be
\omega(e)=\left\{\begin{array}{ll}
   0 & \mbox{ if } h(e)\in\de,\\
   1 & \mbox{ if } h(e)\in\ol\de\setminus\de.
\end{array}\right.
\ee
Let $D$ be the set of edges having both endvertices in $\La_{L+2,M+2}$ which
either are dual to plaquettes in $\de$ or join a vertex of
$\La_{L+1,M+1}$ to a vertex of $\partial\La_{L+1,M+1}$.
We apply Proposition~\ref{Q} to the set $D$, and deduce that
the number of components in the graph $(\ZZ,\eta(\omega))$ having a vertex
in $V(\odelta)$ is simply $K_\delta$. We shall make
use of this observation in the next section when we consider
conditioning on events of the form 
$\{\Delta(\omega)=\de\}$.

\demo{Proof of Proposition~\ref{P}} This may be proved by extending the proof of Lemma~7.2 of 
\cite{GrH}. Instead, we present a variant of that proof.
Consider the set of edges with exactly one endvertex in $V$ and let 
$P$ be the corresponding set of plaquettes.

Let $x\in V$. We first show that $x\in\ins([P])$.
Let  ${\cal U}$ be the set of all closed unit cubes of $\RR$ having centres in 
$V$.   Since all relevant sets in this proof are simplicial, the notions
of path-connectedness and arc-connectedness coincide.
We recall that an unbounded path of $\RR$ from $x$ is 
defined to be a continuous mapping
$\gamma:[0,\infty) \rightarrow \RR$ with $\gamma(0)=x$ whose image is unbounded. 
Any such path $\gamma$ satisfying $|\gamma(t)|\to\infty$
as $t\to\infty$ has a final point $z(\gamma)$
belonging to the (closed) union of all cubes in ${\cal U}$.
Now $z(\gamma)\in[P]$ for all such $\gamma$, and therefore $x\in\ins([P])$.

Let $P_1,P_2,\ldots,P_n$ be the partition of $P$ such that 
the sets $[P_1],[P_2],\ldots,[P_n]$ are the $1$-connected
components of $[P]$ in $\RR$.  Note that $[P_i]\cap [P_j]$ is a finite (or empty) set
for $i\ne j$. 
We show next that there exists $i$ such that $x\in \ins([P_i])$.  Suppose for the sake
of contradiction that this is
false, which is to say that $x\notin \ins([P_i])$ for all $i$. Then
$x\notin\overline{P}_i=[P_i] \cup \ins([P_i])$ for  $i=1,2,\ldots,n$. Note that
each $\ol P_i$ is a closed set which does not separate $\RR$.

Let $i\ne j$. We claim that: either $\overline{P}_i \cap \overline{P}_j$ is a finite set, 
or  one of the sets $\ol P_i$, $\ol P_j$ is a subset of the other. To see this,
suppose that $\ol P_i \cap \ol P_j$ is an infinite set. Suppose further that
$\ol P_i\cap[P_j]$ is infinite. Since $[P_j]$ is a union of unit squares and $\ol P_i$ is a union
of unit squares and cubes, all with corners in $\ZZ+(\frac12,\frac12,\frac12)$, there exists
some edge $f$ of $\Ld_\rd$ such that $f \subseteq \ol P_i\cap[P_j]$. 
We cannot have $f\subseteq [P_i]$ since $[P_i]\cap[ P_j]$ is 
finite, whence $\of \subseteq \ins([P_i])$, where 
$\of$ denotes the open straight-line segment
of $\RR$ joining the endvertices of $f$.
Now $[P_j]$ is 1-connected and $[P_i]\cap[P_j]$ is finite, so that $[P_j]$ is contained in
the closure of $\ins([P_i])$, implying that $[P_j]\subseteq\ol P_i$ and therefore $\ol P_j\subseteq
\ol P_i$.

Suppose next that $\ol P_i\cap[P_j]$ is finite but $\ol P_i\cap\ins([P_j])$ is infinite.
Since $[P_i]$ is 1-connected, it has by definition no finite cutset. Since $[P_i]\cap [P_j]$
is finite, either $[P_i]\subseteq \ol P_j$ or $[P_i]$ is contained in the closure
of the unbounded component of $\RR\setminus [P_j]$. The latter cannot hold since
$\ol P_i\cap\ins([P_j])$ is infinite, whence
$[P_i]\subseteq \ol P_j$ and therefore $\ol P_i \subseteq \ol P_j$.

It follows that we may write $R=\bigcup_{i=1}^{n} \overline{P}_i$ as the union of a collection  of 
closed bounded sets $ \widetilde{P}_i$, 
$i=1,2,\ldots,k$ where $ k\leq n$, that do not separate $\RR$ and  
such that $\widetilde{P}_i\cap  \widetilde{P}_j$  is finite for  $i\neq j$.   
This implies by \cite{Ku} (\S 59, Section II, Theorem 11) that $R$  does not 
separate $\RR$. Now $x\notin R$, whence $x$ lies 
in the unique component of the complement 
$\RR\setminus R$, in contradiction of the assumption that  $x\in\ins([P])$.
We deduce that there exists $k$ such that $x \in\ins([P_k])$, and we define
$Q=P_k$.

Consider now a vertex  $y\in V$. Since $G=(V,E)$ is connected, there exists a path
in $\Ld$ that connects $x$ with $y$ using only vertices in $V$. Whenever $u$ and $v$ are two 
consecutive vertices on this path, $h(\bond{u,v})$  does not belong to $P$. It follows that
$y$ lies in the inside of $[Q]$. 
Claims (i) and (ii) are now proved with $Q$ as given,
and it remains to prove (iii).

Let $W$ be as in (iii), and let $w\in W$.
There exists a path on $\Ld$ from $w$ to infinity using no
vertices of $V$. Whenever $u$ and $v$ are two consecutive vertices on such a path, 
the plaquette $h(\bond{u,v})$ does not lie in $P$. It follows that $w\in\out([P])$,
and therefore $w\in\out([Q])$.
\QED

\demo{Proof of Proposition~\ref{Q}}
Let $H=(\Dv C, \De C)$, and let $H_x$ be the connected component of $H$
containing the vertex $x$. We claim that there exists a plaquette $h_x=h(\bond{y,z})\in\de$
such that $y\in H_x$. This may be shown as follows.

The claim holds with $y=x$ and $h_x=h(\bond{x,z})$
if $x$ has a neighbour $z$ with $h(\bond{x,z})\in\de$. Assume
therefore that $x$ has no such neighbour $z$. There exists a neighbour $u$ of $x$ with
$h(\bond{x,u})\in\odelta\setminus\de$. By a consideration
of the various possibilities,
there exists $\widetilde h\in\de$
such that $\widetilde h\1sim h(\bond{x,u})$, and
\begin{eqnarray*}
\mbox{either} & \mbox{(i)} & \widetilde h=h(\bond{u,z}) \mbox{ for some } z,\\
\mbox{or} & \mbox{(ii)} & \widetilde h=h(\bond{v,z}) \mbox{ for some } v\sim x,\ z\sim u.
\end{eqnarray*}
If (i) holds we take $y=u$, $h_x=\widetilde h$, and if (ii) holds
we take $y=v$ ($\in H_x$), $h_x=\widetilde h$.

We apply Proposition~\ref{P} with $G=H_x$ to obtain a splitting set $Q_x$, and we claim that
\be
\label{qind}
Q_x \cap\de\ne\es.
\ee
This we prove as follows.
If $h_x\in Q_x$, the claim is immediate. Suppose then that
$h_x\notin Q_x$, so that $[h_x]\cap\ins([Q_x])\ne\es$,
implying that $\de$ intersects both $\ins([Q_x])$ and $\out([Q_x])$.
Since $\de$ and $Q_x$ are 1-connected sets of plaquettes, it follows that
$\de\cup Q_x$ is 1-connected. Therefore there exist
$h'\in\de$, $h''\in Q_x$ such that $h'\1sim h''$.
If $h''\in \de$, then (\ref{qind}) holds,
so we may assume that $h''\notin \de$, and hence $h''\in\odelta\setminus\de$.
Then $h''=h(\bond{u,v})$ where $u\in H_x$, and 
therefore $v\in H_x$, a contradiction.
We conclude that (\ref{qind}) holds.

We claim that (\ref{qind}) implies $Q_x\subseteq\de$. Suppose
on the contrary that $Q_x\not\subseteq\de$, so that there
exist $h'\in\de$, $h''\in Q_x\setminus\de$ such that $h'\1sim h''$.
This leads to a contradiction by the argument just given,
whence $Q_x\subseteq\de$.

Suppose now that $x$ and $y$ are vertices of $H$ such that $H_x$ and
$H_y$ are distinct connected components. Then either $H_x$ lies in
$\out([Q_y])$, or $H_y$ lies in $\out([Q_x])$.
Since $Q_x,Q_y\subseteq\de$, either possibility contradicts the
assumption that $x$ and $y$ are connected in $C$.
Therefore $H_x=H_y$ as claimed.
\QED

\demo{Proof of Proposition~\ref{intuition}} 
If $\omega\in{\cal I}_{L,M}\cap\Omega^\mu_{L,M}$, then by definition
$\omega(e)=0$ whenever $h(e)\in\Delta(\omega)$. Suppose conversely that
$\delta\in\DLM$, and let $\omega\in\Omega^\mu_{L,M}$ satisfy $\omega(e)=0$ whenever
$h(e)\in\delta$. Since $\omega \le \ol\omega_\de$, it suffices to show that $\ol\omega_\de\in\ILM$.
Since $\de\in\DLM$, there exists $\xi\in\ILM\cap\Omega^\mu_{L,M}$ such that
$\de=\Delta(\xi)$. Note that $\xi\le \ol\omega_\de$.
Suppose for the sake of a 
contradiction that $\ol\omega_\de\notin\ILM$, and think of
$\ol\omega_\de$ as being obtained from $\xi$ 
by declaring a certain sequence $e_1,e_2,\dots,e_r$
with $\xi(e_i)=0$ for $1\le i \le r$, 
in turn, to be open. Let $\xi^k$ be obtained from $\xi$ by
$\eta(\xi^k)=\eta(\xi)\cup\{e_1,e_2,\dots,e_k\}$. 
By assumption, there exists $K$ such that
$\xi^K\in\ILM$ but $\xi^{K+1}\notin \ILM$. 
For $\psi\in\Omega_{L,M}^\mu$, 
let $J(\psi)$ denote the set of all edges $e$ having endvertices in $\Lambda_{L,M}$, 
with $\psi(e)=1$, and both of whose endvertices
are attainable from $\partial^+\La_{L,M}$ by open paths
of $\psi$.
We apply Proposition~\ref{P} to the finite
connected  graph induced by $J(\xi^K)$
to find that there exists a
splitting set $Q$ of plaquettes such 
that: $\partial^+\La_{L,M}\subseteq \ins([Q])$,
$\partial^-\La_{L,M}\subseteq \out([Q])$,
and $\xi^K(e)=0$ whenever $e\in E_{L,M}$ and $h(e)\in Q$. It must be the case
that $h(e_{K+1})\in Q$, since $\xi^{K+1}\notin \ILM$.
By the 1-connectedness of $Q$, there exists a sequence
$f_1=e_{K+1}, f_2,f_3,\dots,f_t$ 
of edges such that: 
\begin{mylist}
\item $h(f_i) \in Q$ for all $i$,
\item $f_i\in E_{L,M}$ for $1\le i<t$, $f_t=h(\la x,x-(0,0,1)\ra)$
for some $x=(x_1,x_2,1)\in \partial^+\La_{L,M}$,
\item $h(f_i)\1sim h(f_{i+1})$ for $1\le i<t$.
\end{mylist}
It follows that
$h(f_i)\in\de$ for $1\le i\le t$. In particular,
$h(e_{K+1})\in\de$ and so $\ol\omega_\de(e_{K+1})=0$,
a contradiction.  Therefore $\ol\omega_\de\in\ILM$ as claimed.
\QED


\section*{6. Probability distribution of the interface}
For conciseness of notation, we shall henceforth abbreviate $\phi^\mu_{\Lambda_{L,M},p,q}$
to $\phi_{L,M}$, and $\ol\phi^\mu_{\Lambda_{L,M},p,q}$ to $\ol\phi_{L,M}$.
Let $\delta\in\DLM$. We derive next an expression for the probability 
$\phi_{L,M} ( \Delta=\delta)$, which we abbreviate to $\phi_{L,M} (\delta)$.

Let $K_\de$ be the number of components of the graph 
$(\ZZ,\eta(\ol\omega_\de))$, and recall from the discussion after 
Proposition~\ref{Q} that, if $\omega\in\ILM\cap\OLM$ and $\Delta(\omega)=\de$, then
$\omega$ has exactly $K_\de$ open components intersecting $V(\odelta)$.
We have that
\begin{eqnarray}
\label{phidelta}
\displaystyle   
     \phi_{L,M}(\delta)&=&\frac{1}{Z(E_{L,M})} 
           p^{|\overline{\delta}\setminus \delta|}
	   (1-p)^{|\delta|}
     \displaystyle 
       \sum\limits_{{\omega\in\Omega^\mu_{L,M}:}\atop{\Delta(\omega)=\de}}       \biggl\{\prod_{e\in W(\de)}
             p^{\omega(e)}
             (1-p)^{1-\omega(e)}\biggr\}
	     q^{k(\omega)}\nonumber\\
&=&\frac{Z^1(\delta)}{Z(E_{L,M})}p^{|\overline{\delta}\setminus \delta|}
           (1-p)^{|\delta|}q^{K_\de-1},
\end{eqnarray}
where $Z(E_{L,M})=Z^\mu_{\Lambda_{L,M},p,q}$ and $Z^1(\de)=Z^1_{W(\de),p,q}$
as in (\ref{part}).
In this expression and later, for 
$H\subseteq\Hd$, $|H|$ denotes the number of plaquettes in the set 
$H\cap\{h(e): e\in E_{L,M} \}$. 
The term $q^{K_\de-1}$ arises since the application of `1' boundary
conditions to $\de$ has the effect of uniting the boundaries
of the cavities of $\de$, whereby the number of 
clusters diminishes by $K_\de-1$.

We next exploit properties of the partition functions $Z(\cdot)$ 
in order to rewrite 
(\ref{phidelta}). 
For $i=1,2$, let $L_i>0$, $M_i>0$, $\delta_i \in {\cal D}_{L_i, M_i}$, and 
$e_i\in E(\delta_i)\cap E_{L_i,M_i}$, and let
$$
G(e_1, \delta_1, E_{L_1,M_1};e_2, \delta_2, E_{L_2,M_2})
=\sup\Big\{L:\begin{array}[t]{l}
Q_{L}(e_1)\cap E_{L_1,M_1} \simeq Q_{L}(e_2)\cap E_{L_2,M_2} \\
\mbox{and }  Q_{L}(e_1)\cap E(\delta_1) 
	         \simeq Q_{L}(e_2)\cap  E(\delta_2)\Big\},\end{array}
$$
where $Q_L(e)=e+Q_L$ as before. We write $Z^1(E_{L,M})=Z^1_{\Lambda_{L,M},p,q}$.


\begin{pro}
\label{propositionf}
Let $L,M\ge 1$ and $\delta\in \DLM$. We may write $\phi_{L,M}(\de)$ in the form
\be
\label{propositionf1}
     \phi_{L,M}(\delta)
         =\displaystyle\frac{Z^1(E_{L,M})}{Z(E_{L,M})} 
           p^{|\odelta\setminus \delta|}
	   (1-p)^{|\delta|}q^{K_{\delta}-1}
     \exp\left(\sum_{e\in E(\delta)\cap E_{L,M}} f_p(e,\delta,L,M)\right),      
\ee
for functions $f_p(e,\delta, L,M)$ with the following properties.
For $q\ge 1$ there exist $p^\ast<1$ and constants
$C_1$, $C_2$, $\gamma>0$ such that, if $p>p^\ast$,
\begin{eqnarray}
  \label{propositionf2}
       |f_p(e,\delta, L,M)|&<&C_1,\\
\label{propositionf3}
      |f_p(e_1,\delta_1, L_1,M_1)-f_p(e_2,\delta_2, L_2,M_2)| &\leq& C_2 
  e^{-\gamma G},\mbox{ $e_1\in\de_1$, $e_2\in\de_2$, $e_1\simeq e_2$},
\end{eqnarray}
where $G= G(e_1, \delta_1, E_{L_1,M_1};e_2, \delta_2, E_{L_2,M_2})$. 
Inequalities (\ref{propositionf2}) and (\ref{propositionf3}) are 
valid for all relevant values of their arguments.
\end{pro}

\demo{Proof}
We have by Lemma~\ref{logint} that
\be
\label{ratio}
\log\left(\frac{Z^1(\de)}{Z^1(E_{L,M})}\right)
 = \sum_{f\in W(\de)}\Big[g(f,W(\de))-g(f,E_{L,M})\Big]
-\sum_{f\in E(\odelta)}g(f,E_{L,M}),
\ee
where $g(f,D)=g^1_{D,p,q}(f)$.
The summations may be expressed as sums over edges $e$
lying in $E(\delta)$  in the following way.
The edges in $\Ed$ may be ordered according to the lexicographic 
ordering of their
centres. Let $f\in E_{L,M}$ and  $\delta\in \DLM$. Amongst
all edges in $E(\delta)\cap E_{L,M}$ which are closest to $f$ (in the sense
that their centres are closest in $L^\infty$ norm),  we write
$\nu(f,\de)$ for the earliest edge in this ordering.
We have  by (\ref{ratio}) that
\be
\label{ratio2}
\log\left(\frac{Z^1(\de)}{Z^1(E_{L,M})}\right)
 = \sum_{e\in E(\de)\cap E_{L,M}}
 f_p(e,\delta,L,M)
\ee
where
\be
\label{defoff}
f_p(e,\delta,L,M)=\displaystyle
\sum\limits_{{f\in W(\de):}\atop{\nu(f,\de)=e}}
  \Big[g(f,W(\de))-g(f,E_{L,M})\Big] -
\sum\limits_{{f\in E(\odelta):}\atop{\nu(f,\de)=e}}
  g(f,E_{L,M}).
\ee
This establishes (\ref{propositionf1}) via (\ref{phidelta}).

It remains to show the required properties of the $f_p$.
Suppose $e=\nu(f,\de)$ and set $r=\Vert e,f\Vert$. Then $Q_{r-2,r-2}(f)$ does
not intersect $\odelta$, implying by Lemma~\ref{expbnd} that
\be
\label{bnd1}
|g(f,W(\de))-g(f,E_{L,M})|\le e^{-\alpha\Vert e,f\|+2\alpha}
\quad\mbox{if } p>p^\ast,
\ee
where $p^\ast$ and $\alpha$ are given as in that lemma.
Secondly, there exists an absolute constant $K$ such that,
for all $e$ and $\de$, the number
of edges $f\in E(\odelta)$ with $e=\nu(f,\de)$ is no greater than $K$. Therefore,
by (\ref{bndg}),
$$
   |f_p(e,\delta, L,M)|
       \leq\sum_{f\in \Ed}e^{-\alpha\|e,f\|+2\alpha}
                + K  (1-p)(q-1) 
$$
as required for (\ref{propositionf2}).

Finally we show  (\ref{propositionf3}) for $p>p^*$ and appropriate $C_2$, $\gamma$.
Let $e\in\de_1$, $e_2\in\de_2$, and let
$G$ be given as in the proposition; we may 
suppose that $G>9$. By assumption, $e_1\simeq e_2$, whence there exists a translate
$\tau$ of $\Ld$ such that $\tau e_1=e_2$. We have for $f\in W(\de_1)\cap Q_{G/3}(e_1)$ that
\begin{eqnarray}
\label{bnd2}
\tau\Big[Q_{G/3}(f)\cap E_{L_1,M_1}\Big]&=&
  Q_{G/3}(\tau f)\cap E_{L_2,M_2},\\
\label{bnd3}
\tau\Big[Q_{G/3}(f)\cap \de_1\Big]&=&
  Q_{G/3}(\tau f)\cap \de_2,
\end{eqnarray}
and that
\be
\label{bnd4}
\mbox{ for } \| f,e_1\| \le \tfrac13 G,\ 
\nu(f,\delta_1)=e_1 \mbox{ if and only if } \nu(\tau f,\delta_2)=e_2.
\ee

It follows from the definition (\ref{defoff}) of the functions $f_p$ that
\begin{eqnarray}
 \label{bigsum}
 &&|f_p(e_1,\delta_1, L_1,M_1)-f_p(e_2,\delta_2, L_2,M_2)| \nonumber \\  
 &&\displaystyle\hspace{.5cm}\leq \sum_{{f\in W(\de_1)\cap Q_{G/3}(e_1):}\atop
                               {\nu(f,\delta_1)=e_1}}
                \Big\{|g(f,W(\de_1))- g(\tau f,W(\de_2))|+
                      |g(f,E_{L_1,M_1})- g(\tau f,E_{L_2,M_2})|\Big\}\nonumber \\	
      &&\displaystyle\hspace{3cm}+\sum_{{f\in W(\de_1)\setminus Q_{G/3}(e_1):}\atop
                               { \nu(f, \delta_1)=e_1}}
                |g(f,W(\de_1))- g(f,E_{L_1,M_1})|\nonumber \\
      &&\displaystyle\hspace{3cm}+\sum_{{f\in W(\de_2)\setminus Q_{G/3}(e_2):}\atop
                                { \nu(f,\delta_2)=e_2}}
                |g(f,W(\de_2))- g(f,E_{L_2,M_2})|+S,
\end{eqnarray}
where
$$
S=\left\vert \sum_{{f\in E(\odelta_1):}\atop{\nu(f,\delta_1)=e_1 }}
                         g(f,E_{L_1,M_1})-\sum_{{f\in E(\odelta_2):}\atop{\nu(f,\delta_2)=e_2 }}
g(f,E_{L_2,M_2})\right|.
$$

By (\ref{bnd2}), (\ref{bnd3}), and  Lemma~\ref{expbnd}, the first summation 
in (\ref{bigsum}) is bounded
above by $2G^3e^{-\frac13 \alpha G}$.
Using the definition of the $\nu(f,\de_i)$,
the second and third summations of (\ref{bigsum}) are bounded above, respectively, by
$$
\sum_{f\notin Q_{G/3}(e_i)} e^{-\alpha\|f,e_i\|+2\alpha}
\le C'e^{-\frac13\alpha G +2\alpha},
$$
for some $C'<\infty$, as in (\ref{bnd1}).
We have by (\ref{bnd4}) that
$$
S= \left\vert \sum_{{f\in E(\odelta_1):}\atop{\nu(f,\delta_1)=e_1 }}
                         g(f,E_{L_1,M_1})-
g(\tau f,E_{L_2,M_2})\right|
\le Ke^{-\frac13\alpha G},
$$
and inequality (\ref{propositionf3}) is proved for an appropriate choice of $\gamma$.
\QED

In the next part of this section, we consider measures and interfaces for
the infinite cylinder $\La_L=\La_{L,\infty}=[-L,L]^2\times\Z$. We note first that,
if $q\ge 1$, then $\phi_{L,M+1}\lest \phi_{L,M}$, as in \cite{Gr1}, Theorem 3.1(a),
whence the (decreasing) weak limit
\be
\label{lim exists}
\phi_L = \lim_{M\to\infty}\phi_{L,M}
\ee
exists. We write $\Omega_L^\mu$ for the set of all configurations $\omega$ such that
$\omega(e)=\mu(e)$ for $e\notin E_L=\lim_{M\to\infty}E_{L,M}$, and $\IL$
for the event that no vertex of $\partial\La_L^+$ is joined by an open path to a vertex
of $\partial\La_L^-$.  The set of interfaces on which we concentrate is
$\DL=\bigcup_M\DLM=\lim_{M\to\infty}\DLM$. Thus $\DL$ is 
the set of interfaces which span
$\La_L$, and every member of $\DL$ is bounded in the direction of the third
coordinate. It is easy to see that
$\IL\supseteq\lim_{M\to\infty}\ILM$,
and it is a consequence of the next
lemma that the difference between these two events has $\phi_L$-probability zero.


\begin{lem}
\label{limitphi}
We have, if $q \ge 1$, that $\phi_{L,M}(\cdot\mid\ILM)\Rightarrow
\phi_L(\cdot\mid\IL)$ as $M\to\infty$, and that
$$
\phi_L\Big(\IL\setminus\lim_{M\to\infty}\ILM\Big)=0. 
$$
\end{lem}

For $L_i>0$, $\delta_i\in {\cal D}_{L_i}$, and $e_i\in E(\delta_i)\cap E_{L_i}$, 
let
$$
G(e_1, \delta_1, E_{L_1};e_2, \delta_2, E_{L_2})
      =\sup\Big\{L: \begin{array}[t]{l}
    Q_{L}(e_1)\cap E_{L_1} \simeq Q_{L}(e_2)\cap E_{L_2} \\
        \mbox{and }  Q_{L}(e_1)\cap E(\delta_1)\simeq Q_{L}(e_2)\cap  E(\delta_2)\Big\}.
\end{array}
$$
On the event $\IL$, $\Delta$ is defined as before to be the maximal
1-connected set of open plaquettes which intersects $\de_0\setminus E_L$.

\begin{lem}
\label{limitf}
\mbox{\rm(a)} Suppose $L>0$, $\delta\in \DL$, and $e\in E(\delta)\cap E_L$. 
The functions $f_p$ given in (\ref{defoff}) are such that the limit
\be
\label{flimit}
  f_p(e,\delta,L)=\lim_{M\to \infty} f_p(e,\delta,L,M)
\ee
exists. Furthermore, if $p>p^*$,
\be
\label{cbound}
 |f_p(e,\delta, L)|<C_1, 
\ee
and, for $L_i>0$, $\delta_i\in {\cal D}_{L_i}$, 
and $e_i\in E(\delta_i)\cap E_{L_i}$ satisfying $e_1\simeq e_2$, 
$$
      |f_p(e_1,\delta_1, L_1)-f_p(e_2,\delta_2, L_2)|
       \leq C_2 e^{-\gamma G},
$$
where $p^*$, $C_1$, $C_2$, $\gamma$ are given as in Proposition~\ref{propositionf}
and $G=G(e_1, \delta_1, E_{L_1};e_2, \delta_2, E_{L_2})$.
\par
\noindent
\mbox{\rm(b)} For $q\ge 1$ and $\de\in\DL$, the probability
$\phi_L(\de\mid\IL)=\phi_L(\Delta=\de\mid\IL)$ is given by
\be
\label{probdelta}
\phi_L(\de\mid\IL)=\frac 1{Z_L} p^{|\odelta\setminus\de|}(1-p)^{|\de|}
q^{K_\de}\exp\left(\sum_{e\in E(\de)\cap E_L} f_p(e,\de,L)\right),
\ee
where $Z_L$ is the appropriate normalizing constant.
\end{lem}

\demo{Proof of Lemma~\ref{limitphi}} 
It suffices for the claim of weak convergence that
\be
\label{suff}
\phi_{L,M}(F\cap\ILM)\to \phi_L(F\cap\IL)\quad\mbox{ for all cylinder events } F.
\ee
Let $A_{L,M}=[-L,L]^2\times\{-M\}$ and $B_{L,M}=[-L,L]^2\times\{M\}$, and let
$T_{L,M}$ be the event that no open path exists between a vertex of
$\partial\La_{L,M}^+\setminus B_{L,M}$ and  a vertex of
$\partial\La_{L,M}^-\setminus A_{L,M}$. Note that $T_{L,M}\to\IL$
as $M\to\infty$.
Let $F$ be a cylinder event. Then
\begin{eqnarray}
\label{uplim}
\phi_{L,M}(F\cap\ILM) &\le& \phi_{L,M}(F\cap T_{L,M'}) \quad\mbox{for } M'\le M\nonumber\\
 &\to& \phi_L(F\cap T_{L,M'}) \quad\quad\mbox{as } M\to\infty\nonumber\\
&\to& \phi_L(F\cap \IL) \quad\quad\quad\hskip1mm\mbox{as } M'\to\infty.
\end{eqnarray}

In order to obtain a corresponding lower bound, we introduce the event
$K_r$ that all edges of $E_L$, both of whose endvertices
have third coordinate equal to $\pm r$, are open. We may suppose without
loss of generality that $p>0$.
We have by Lemma \ref{comp} 
that $\phi_{L,M}$ dominates product measure with density $\pi=p/\{p+(1-p)q\}$,
whence there exists $\beta=\beta_L<1$ such that 
$$
\phi_{L,M}(K_r\mbox{ for some } r\le R) \ge 1-\beta^R\q\mbox{for } R<M.
$$
Now $\ILM\subseteq T_{L,M}$, and $T_{L,M}\setminus \ILM
\subseteq \bigcap_{r=1}^{M-1} K_r^\rc$, whence
\begin{eqnarray}
\label{downlim}
\phi_{L,M}(F\cap \ILM) &\ge& \phi_{L,M}(F \cap T_{L,M}) - \beta^{M-1}\nonumber\\
&\ge& \phi_{L,M}(F\cap \IL) - \beta^{M-1}\nonumber\\
&\to& \phi_L(F\cap \IL)\quad\quad\quad\mbox{as } M\to\infty.
\end{eqnarray}
Equation (\ref{suff}) follows from (\ref{uplim}) and (\ref{downlim}).
The second claim of the lemma follows by taking $F=\Omega$,
the entire sample space.
\QED

\demo{Proof of Lemma~\ref{limitf}} (a) The existence of the 
limit follows from the monotonicity
of $g(f,D_i)$ for an increasing sequence $\{D_i\}$, and the proof of
(\ref{propositionf2}). The inequalities are implied by 
(\ref{propositionf2}) and (\ref{propositionf3}).

\noindent
(b) Let $\de\in\DL$, so that 
$\de\in\ILM$ for all large $M$.
By Lemma \ref{limitphi}, $\phi_L(\de\mid\IL) = 
\lim_{M\to\infty}\phi_{L,M}(\de\mid\ILM)$.  We take the limit
as $M\to\infty$ in (\ref{propositionf1}), and use part (a) to obtain the claim.
\QED


\section*{7. Geometry of interfaces}
Next, we describe in more detail the interfaces in $\DL=\lim_{M\to\infty}\DLM$. 
While it was natural in Section 5 to
introduce the extended interface $\odelta$ of a member $\delta$ of $\DL$, it turns
out to be useful when studying its geometry to introduce its
{\it semi-extended interface\/}  
$$
    \sdelta= \delta\cup \Big\{h\in \Hd: h \smbox{is a horizontal plaquette that is
                     $1$-connected to}\delta \Big\}.
$$
Let $x=(x_1,x_2,x_3)\in\ZZ$.
The {\it projection\/} $\pi(h)$ of a
horizontal plaquette $h=h(\bond{x,x+(0,0,1)})$ onto the regular interface
$\delta_0$ is defined to be the plaquette 
$\pi(h)= h(\bond{(x_1,x_2,0),(x_1,x_2,1)})\in \delta_0$. The projection of the vertical plaquette
$h=h(\bond{x, x+(1,0,0)})$ is the interval 
$\pi(h)=[(x_1+\frac12 , x_2-\frac12 , \frac12 ),
(x_1+\frac12 , x_2+\frac12 , \frac12 )]$. Similarly,   
$h=h(\bond{x, x+(0,1,0)})$ has projection
$\pi(h)=[(x_1-\frac12 , x_2+\frac12 , \frac12 ),(x_1+\frac12 , x_2+\frac12 ,
\frac12 )]$.

Let $\de\in\DL$.
A horizontal plaquette $h$ of the 
semi-extended interface $\sdelta$ is called a {\it c-plaquette\/} if 
$h$ is the unique member of $\sdelta$ having projection $\pi(h)$.
All other plaquettes of  $\sdelta$ are called {\it w-plaquettes}.
A {\it ceiling\/} of $\de$ is a maximal $0$-connected set of c-plaquettes. The 
{\it projection\/} of a ceiling
$C$ is the set 
$\pi(C)=\{\pi(h): h \in C\}$.
Similarly, we define  a {\it wall\/} $W$ of $\delta$ as  a maximal $0$-connected set of 
w-plaquettes, and its {\it projection\/} as
$$
   \pi(W)=\{\pi(h): h \smbox{is a horizontal plaquette of} W\}.
$$ 
We collect together some properties of interfaces thus.


\begin{lem}
\label{properties}
Let $\delta\in \DL$.
\begin{mylist}
\item The set $\sdelta\setminus\de$ contains no  c-plaquette.
\item
All plaquettes of $\sdelta$ that are $1$-connected to some c-plaquette
 are horizontal plaquettes of 
$\delta$. All horizontal plaquettes that are $0$-connected to some c-plaquette
belong to $\sdelta$. 
\item Let $C$ be a ceiling. There is a unique plane parallel to the regular
interface which contains all the c-plaquettes of $C$.
The set of all horizontal plaquettes, which are $0$-connected to
members of $C$ but do not themselves lie in $C$, form
a $0$-connected subset of $\sdelta$.
\item Let $C$ be a ceiling.  We have that
   $\{h\in \sdelta: \pi(h)\subseteq [\pi(C)]\}=C$.
\item Let $W$ be a wall. We have that
   $\{h\in \sdelta: \pi(h)\subseteq [\pi(W)]\}=W$.
\item For each wall $W$, 
$\delta_0\setminus \pi(W)$ has exactly one maximal infinite 
$0$-connected component (respectively, $1$-connected component). 
\item Let $W$ be a wall, and 
suppose that $\delta_0\setminus  \pi(W)$ comprises $n$ maximal 
$0$-connected sets $H_1,H_2, \ldots,H_n $.  The set of all plaquettes 
$h\in \sdelta\setminus W$ 
which are $0$-connected to $W$ comprises only c-plaquettes,
which belong to the union of exactly $n$ distinct
ceilings  $C_1,C_2\ldots,C_n$ such that 
$\{\pi(h):\mbox{\rm\ $h$ is a c-plaquette of $C_i$} \}\subseteq H_i$.
\item The projections $\pi(W_1)$ and $\pi(W_2)$ of two different walls 
$W_1$ and $W_2$ of  $\sdelta$ are not $0$-connected.
\item The projection $\pi(W)$ of any wall $W$ contains at least one plaquette of $\delta_0$. 
\end{mylist}
\end{lem}

The displacement of the plane in (iii) from
the regular interface, counted positive or negative, 
is called the {\it height\/} of the ceiling $C$.

\demo{Proof}
(i) Let $h$ be a c-plaquette of $\sdelta$ with $\pi(h)=h_0$.   
Since $\delta\in\DL$, it contains at least one plaquette with projection $h_0$.
Yet, according to the definition of a c-plaquette, there is no such a 
plaquette besides $h$.  Therefore $h\in\de$.
  
\noindent
(ii) Suppose $h$ is a c-plaquette. Necessarily, $h$ belongs to  $\delta$ and any horizontal 
plaquette which is $1$-connected to $h$ belongs to $\sdelta$.
It may be seen in addition that any vertical plaquette which is
1-connected to $h$ lies in $\ol\de\setminus\de$.
Suppose, on the contrary, that one such vertical plaquette $h'$ lies
in $\de$. Then the horizontal plaquettes which are 1-connected
to $h'$ lie in $\sdelta$. One of these latter
plaquettes has projection $\pi(h)$, in contradiction of the assumption that
$h$ is a  c-plaquette. 

We may now see as follows that any horizontal plaquette $h''$ which
is $1$-connected to $h$ must lie in $\de$.
Suppose, on the contrary, that one such plaquette $h''$ lies in $\ol\de\setminus\de$.
We may construct a path of open edges on $(\ZZ, \eta(\ul\omega_\delta ))$ that connects the vertex 
$x$ just above $h$ with the vertex $x-(0,0,1)$ just below $h$, using the open edges 
of $\ul\omega_\de$ corresponding to the three relevant
plaquettes of $\ol\de\setminus\de$. This contradicts the assumption that
$h$ is a c-plaquette of the interface $\de$.

The second claim of (ii) follows immediately, by the definition of $\sdelta$.

\noindent(iii)
The first part follows by the 
definition of ceiling, since the only horizontal plaquettes that are 
$0$-connected with a given c-plaquette $h$  lie in the plane containing $h$.
The second assertion follows from (ii) and the geometry of $\Z^2$.

\noindent(iv)
Assume that $h\in\sdelta$ and $\pi(h)\subseteq[\pi(C)]$.
If $h$ is horizontal, the conclusion holds by the definition of c-plaquette.
If $h$ is vertical, then $h\in\de$,  and all 1-connected horizontal plaquettes lie in
$\sdelta$. At least two such horizontal plaquettes project onto the same plaquette 
in $\pi(C)$, in contradiction of the assumption that $C$ is a ceiling. 

\noindent(v)
Let $C$ be a ceiling and let $\gamma_1, \gamma_2,\ldots, \gamma_n$ 
be the maximal  $0$-connected sets of plaquettes  of
$\delta_0\setminus \pi(C)$.
Let 
$\delta_i^\ast=\{h\in \sdelta: \pi(h)\subseteq [\gamma_i] \}$.
We have by part (iv) that
 $\sdelta=C\cup (\bigcup_{i=1}^n \delta_i^\ast)$. We claim that 
each $\delta^\ast_i$ is $0$-connected, and we prove this as follows.  
Let $h_1,h_2\in \delta^\ast_i$. Since  $\sdelta$ is $0$-connected, 
it contains a 
sequence $h_1=f_0, f_1,\ldots,f_m=h_2$ of plaquettes 
such that 
$f_{i-1}\0sim f_i$  for $1\le i\le m$. We need to show that such
a sequence exists containing no plaquettes in $C$. Suppose
on the contrary that the sequence $(f_i)$ has a non-empty
intersection  with $C$.
Let $k=\min\{i:f_i\in C\}$ and $l=\max\{i:f_i\in C\}$, and
note that
$0<k\leq l<n$.

If $f_{k-1}$ and $f_{l+1}$ are horizontal, by (ii) and (iii),
they are 0-connected by a path of horizontal plaquettes
of $\sdelta\setminus C$, in which case $h_1$ and $h_2$ 
are in the same $0$-connected component of
$\sdelta\setminus C$. A similar argument is valid if either or both
of  $f_{k-1}$ and $f_{l+1}$ is vertical. For example,
if $f_{k-1}$ is vertical, by (ii) it cannot be $1$-connected to 
a plaquette of $C$. Hence it is $1$-connected to some 
horizontal plaquette  in $\sdelta\setminus C$ which is itself 
$1$-connected to a plaquette of $C$. The same conclusion is
valid for  $f_{l+1}$ if vertical. In any such case, 
by (ii) and (iii) 
there exists a $0$-connected  sequence of 
w-plaquettes connecting $f_{k-1}$ with 
$f_{l+1}$, and the claim follows.

To prove (v), we note by the above that the wall $W$ is a
subset of one of the sets $\delta_i^\ast$, 
say $\delta_1^\ast$. Next
we let $C_1$ be a ceiling contained in $\delta_1^\ast$, if this exists, 
and we repeat the above procedure.
We consider the $0$-connected components of $\gamma_1\setminus \pi(C_1)$, 
and we use the fact that $\delta_1^\ast$ is $0$-connected to 
deduce that the set of plaquettes which project  
onto one of these components is itself $0$-connected. 

This procedure is repeated until all ceilings have been removed,
the result being a $0$-connected set of w-plaquettes 
of which, by definition of a wall, all members belong to $W$. 

Finally, (vi) is a simple observation since walls are finite.
 Claim (vii) is immediate from claim (ii) and the definitions of wall and ceiling.
Claim (viii) follows from (v) and (vii), and (ix) is a consequence
of the definition of the semi-extended  interface $\sdelta$. 
\QED

The properties described in Lemma~\ref{properties} allow us to  describe a wall $W$ in more detail. 
By (vi) and (vii), there exists a unique ceiling that is 
$0$-connected to $W$ and with projection in
the infinite 0-connected component of $\delta_0\setminus \pi(W)$. 
We call this ceiling the {\it base\/} of $W$. 
The {\it altitude\/} of $W$ is the height of the base of $W$; see (iii).
The {\it height\/} $D(W)$ of $W$ is the maximum 
absolute value of the displacement
in the third coordinate direction of $[W]$ from the
horizontal plane $\{(x_1,x_2,s+\tfrac12): x_1,x_2\in \Z\}$, where $s$ is
the altitude of $W$. 
The {\it interior\/} 
$\inter(W)$ (of the projection $\pi(W)$) of $W$ is the 
complement in $\de_0$ of the unique maximal infinite 0-connected component
of $\de_0\setminus\pi(W)$ (cf.\ (vi)).

We next define the concept of a standard wall. Let $S=(A,B)$
where $A$, $B$ are sets of plaquettes. We call $S$
a {\it standard wall\/} if there exists
$\de\in\DL$ such that $A\subseteq\de$, $B\subseteq\sdelta\setminus\de$, and
$A\cup B$ is the unique
wall of $\de$. If $S=(A,B)$ is a standard wall,
we shall refer to plaquettes of either $A$ or $B$ as plaquettes
of $S$, and we write $\pi(S)=\pi(A\cup B)$.

\begin{lem}
\label{unique}
Let $S=(A,B)$ be a standard wall. There exists a unique
$\de\in\DL$ such that\/\mbox{\rm:} $A\subseteq\de$, $B\subseteq\sdelta\setminus\de$, and
$A\cup B$ is the unique
wall of $\de$. 
\end{lem}

This will be proved soon. We denote by $\de_S$ the unique such $\de\in\DL$ corresponding
to the standard wall $S$.
We shall see that
standard walls are the basic building blocks for a general interface. 
Notice that the base of a
standard wall is a subset of the regular interface.  
We introduce an ordering on the plaquettes of $\delta_0$,
and we define the {\it origin\/} of the standard wall $S$
to be the earliest plaquette 
in $\pi(S)$ which is $1$-connected to some plaquette of 
$\delta_0\setminus \pi(S)$.
Such an origin exists by Lemma~\ref{properties}(ix), and the origin
belongs to $S$ by (ii). For $h\in \delta_0$, we denote by 
${\cal S}_h$ the set of all standard walls with origin $h$. We 
attach to ${\cal S}_h$ the {\it empty wall\/} ${\cal E}_h$
interpreted as a wall with origin $h$
but containing no plaquettes.

A family $\{S_i=(A_i,B_i): 1\le i\le m\}$ of standard walls is called  
{\it admissible\/} if:
\begin{mylist}
\item for $i\neq j$, there exists no pair  $h_1\in \pi(S_i)$ and  
$h_2\in \pi(S_j)$
such that  $h_1\0sim h_2$, 
\item if, for some $i$,  
$h(e)\in S_i$ where $e\notin E_L$, then $h(e)\in A_i$ if and only if
$\mu(e)=0$.
 \end{mylist}
The members of any such family have distinct origins. For our future convenience
we label each $S_i$ according to its origin $h(i)$, and write $\{S_h:h\in\de_0\}$
for the family, where $S_h$ is to be interpreted as ${\cal E}_h$ when $h$
is the origin of none of the $S_i$.
We adopt the convention that, when a standard wall is denoted as
$S_h$ for some $h\in\de_0$, then $S_h\in{\cal S}_h$.

We introduce next the concept of a group of walls.
Let $h\in \delta_0$, $\delta\in \DL$, and denote by $\rho(h, \delta)$
the number of (vertical or horizontal) plaquettes in $\de$ whose
projection is a subset of $h$.
Two standard walls $S_1$, $S_2$ are called 
{\it close\/} if there exist
$h_1\in \pi(S_1)$ and $h_2\in \pi(S_2)$ such that
$$
\|h_1,h_2\|<\sqrt{\rho(h_1, \delta_{S_1})}+ \sqrt{\rho(h_2, \delta_{S_2})}.  
$$
A family $G$ of non-empty standard walls 
is called a {\it group of (standard) walls\/} if it is
admissible and if, for
any pair $S_1, S_2 \in G$, there exists a 
sequence $T_0=S_1,T_1, T_2,\ldots,T_n=S_2$ of
members of $G$ such that
$T_i$ and $T_{i+1}$ are close for $0\le i <n$.

The {\it origin\/} of a group of walls
is defined as the earliest of the origins of the standard walls
therein. We write ${\cal G}_h$ for the set
of all possible groups of walls with origin $h\in \delta_0$. As before,
we attach to ${\cal G}_h$ the {\it empty group\/} 
with origin $h$ but containing no standard wall which we denote
also as ${\cal E}_h$.
A family  $\{G_i:1\le i\le m\}$ of groups of walls  is
called {\it admissible\/} if, for $i\ne j$,  there exists no pair 
$S_1\in G_i$, $S_2\in G_j$ such that $S_1$ and $S_2$
are close.

We adopt the convention that, when a group of walls is denoted as
$G_h$ for some $h\in\de_0$, then $G_h\in{\cal G}_h$. Thus
a family of groups of walls may be written as a collection $\GG=\{G_h:h\in\de_0\}$
where $G_h\in{\cal G}_h$.


\begin{lem}
\label{construction}
The set $\DL$  is in one--one correspondence with both the collection of 
admissible families of standard walls, and with the
collection of admissible
families of groups of walls.
\end{lem}

Equally important to the existence of these one--one correspondences
is their nature, as described in the proof of the lemma.
We write $\de_G$ (respectively $\de_{\GG}$) for the interface corresponding
thus to an admissible family $G$ of standard walls (respectively an admissible
family $\GG$ of groups of walls).

\demo{Proof of Lemma~\ref{unique}}
Let $\de\in\DL$ have unique wall $S=(A,B)$. By definition, every plaquette
of $\sdelta$ other than those in $A \cup B$ is a c-plaquette, so that
$\Sigma=\sdelta\setminus(A\cup B)$ is a union of ceilings
$C_1,C_2,\dots,C_n$. Each $C_i$ contains 
some plaquette $h_i$ which is 1-connected
to some $h_i'\in A$, whence, by Lemma~\ref{properties}(iii),
the height of $C_i$ is determined uniquely by knowledge
of $S$.
Hence $\de$ is unique.
\QED

\demo{Proof of Lemma~\ref{construction}}
Let $\delta\in \DL$.
Let $W_1,W_2, \ldots, W_n$ be the non-empty walls of $\sdelta$, and 
write $W_i=(A_i,B_i)$ where $A_i=W_i\cap\delta$, 
$B_i=W_i\cap(\sdelta\setminus\de)$.
Let $s_i$
be the altitude of $W_i$.
We claim that $\tau_{(0,0,-s_i)}W_i$
is a standard wall, and we prove this as follows. 
Let $C_{i_j}$, $j=1,2,\ldots,k$, be the ceilings that are $0$-connected 
to $W_i$, and 
let $H_{i_j}$ be the maximal $0$-connected set of 
plaquettes in $\delta_0\setminus \pi(W_{i})$ onto which 
$C_{i_j}$ projects. (See Lemma~\ref{properties}(vii).)  
It suffices to construct an
interface $\delta({W_i})$  having $\tau_{(0,0,-s_i)}W_i$ as its unique wall. 
To this end
we add to $\tau_{(0,0,-s_i)}A_i$ the plaquettes in $\tau_{(0,0,-s_i)}C_{i_j}$, 
$j=1,2,\ldots,k$, together with, for each $j$, 
the horizontal plaquettes in the maximal $0$-connected set of horizontal plaquettes 
that contains $\tau_{(0,0,-s_i)}C_{i_j}$ and 
elements of which project onto $H_{i_j}$. 

We now define the family $\{S_h: h\in \delta_0\}$ of standard walls 
by
$$
  S_h=\left\{\begin{array}{ll}
\tau_{(0,0,-s_i)} W_i & \smbox{if $h$ is the origin of $\tau_{(0,0,-s_i)}W_i$,}\\
{\cal E}_h  & \smbox{if $h$ is the origin  of no $\tau_{(0,0,-s_i)}W_i$.}
\end{array}\right.
$$ 
More precisely, in the first case,
$S_h=(A_h,B_h)$ where $A_h=\tau_{(0,0,-s_i)}A_i$
and $B_h=\tau_{(0,0,-s_i)}B_i$.
That this is an admissible family of  standard walls 
follows from Lemma \ref{properties}(viii)
and from the observation that $s_i=0$ when $E(W_i)\cap E_L^\rc\neq \es$. 

Conversely, let  $\{S_h=(A_h,B_h):h\in\de_0\}$ be an admissible family
of standard walls.
We shall show that there is a unique  
interface $\delta $ 
corresponding in a certain way to this family.
Let $S_1,S_2 \ldots, S_n$ be the non-empty walls of the family, and let  
$\delta_i$ be the unique 
interface in $\DL$ having $S_i$ as its only wall.

We introduce the partial ordering on the walls given by  $S_i<S_j$ if 
$\inter(S_i)\subseteq \inter(S_j)$, and we re-order the non-empty walls in such a way
that $S_i<S_j$ implies $i<j$. 

When  it exists, we take the first index $k>1$ such that
$S_1<S_k$ and we modify 
$\delta_k$ as follows.
First we  remove the c-plaquettes that  project onto  $\inter(S_1)$,
and then we add translates of the plaquettes of $A_1$.
This is done by translating these
plaquettes so that the base of $S_1$ is raised (or lowered) to the  plane containing 
the ceiling that is $0$-connected 
to   $S_k$ and that projects on the  maximal $0$-connected set  of plaquettes in 
$\delta_0\setminus \pi(S_k)$ that contains   $\pi(S_1)$. 
(See Lemma~\ref{properties}(viii).) We write $\de_k'$ for the ensuing interface.
We now repeat this procedure starting  from the set of 
standard walls $S_2,S_3, \ldots, S_n$ 
and interfaces  $\delta_2,\delta_3, 
\ldots,\de_{k-1},\de_k',\de_{k+1},\ldots,\delta_n$.
If no such $k$ exists, we continue the procedure with
the interfaces $\delta_2,\delta_3, 
\ldots,\de_{k-1},\de_k,\de_{k+1},\ldots,\delta_n$.

We continue this process until we are left with interfaces 
$\delta_{i_k}''$, $k=1,2,\ldots,r$, having indices 
which refer to standard walls that are smaller 
than no other wall. The final interface $\delta$ is now constructed as 
follows. For each $k$, we remove from the regular interface  
$\delta_0$  all horizontal plaquettes contained in $\inter(S_{i_k})$, and we replace  
them by the plaquettes of $\delta_{i_k}''$ that project onto $\inter(S_{i_k})$. 

The final assertion concerning admissible families of groups
of walls is straightforward.
\QED

Next we derive certain 
combinatorial properties of walls.
For $S=(A,B)$ a standard wall, we write $N(S)=|A|$ and we
set
$\Pi(S)=N(S)-|\pi(S)|$.
For an admissible set $F=\{S_1,S_2,\ldots,S_m\}$
of  standard walls, we write
$\Pi(F)=\sum_{i=1}^{m}\Pi(S_i)$, $N(F)=\sum_{i=1}^{m}N(S_i)$,
and $\pi(F)=\bigcup_{i=1}^{m}\pi(S_i)$.

\begin{lem}
\label{Pi}
Let $S=(A,B)$ be a standard wall, and $D(S)$ its height.
\begin{mylist}
\item    $N(S)\geq \frac{14}{13}|\pi(S)|$. Consequently,
$\Pi(S)\geq \frac{1}{13}|\pi(S)|$ and $\Pi(S)\geq \frac{1}{14}N(S)$.
\item $N(S)\geq \frac{1}{5}|S|$.
\item $\Pi(S)\ge D(S)$.
\end{mylist}
\end{lem}
    
\demo{Proof}
(i) Define for each $h_0\in \delta_0$ the set
$ U(h_0)=\{h\in \delta_0: h=h_0 \smbox{or} h\1sim h_0 \}$.
We call two plaquettes $h_1, h_2\in \delta_0$ {\it separated} if $U(h_1)\cap U(h_2)=  \es$.
Denote by $H_{\mbox{\rm \scriptsize sep}}=
H_{\mbox{\rm \scriptsize sep}}(S)\subseteq \pi(S)$ a
set of pairwise-separated plaquettes in $\pi(S)$ having maximum cardinality,
and let $H=\bigcup_{h_1\in \hsep} [U(h_1)\cap\pi(S)]$.
Note that
\be
\label{hbnd}
|\hsep | \ge \tfrac1{13}|\pi(S)|.
\ee

For every $h_0\in \pi(S)$, there exists a horizontal plaquette $h_1\in \de_S$ such that
$\pi(h_1)=h_0$. Since $A\cup B$ contains no
c-plaquette of $\de_S$, it is the case that $h_1$ is
a w-plaquette, whence $h_1\in A$. In particular, $N(S)\geq |\pi(S)|$.

For $h_0=\pi(h_1)\in H_{\mbox{\rm \scriptsize sep}}$ where $h_1\in A$,
we claim that
\be
\label{Piproof1}
\Big|\Big\{h\in A: \mbox{either } \pi(h)\subseteq [h_0]
   \smbox{or} \pi(h)\in U(h_0)\Big\}\Big|
      \geq |U(h_0)\cap \pi(S)|+1.
\ee
It follows from (\ref{hbnd}) and (\ref{Piproof1}) that
\begin{eqnarray*}
N(S) &\ge& \sum_{h_0\in\hsep}\Big\{|U(h_0)\cap \pi(S)| +1\Big\} + |\pi(S)\setminus H|\\
&=& |H| + |\hsep| + |\pi(S)| - |H|
\ge \tfrac{14}{13}|\pi(S)|
\end{eqnarray*}
as required.

In order to prove (\ref{Piproof1}), we argue
first that $U(h_0)\cap \pi(S)$ contains at least one (horizontal) plaquette besides
$h_0$. Suppose that this is not true.
Then $U(h_0)\setminus h_0$ contains the projections of c-plaquettes of $\delta^\ast_S$
only.
By Lemma~\ref{properties}(ii, iii),
these c-plaquettes belong to the same ceiling $C$
and therefore lie in the same plane.
Since $h_1$ is by assumption
a w-plaquette, there must be at least one other horizontal plaquette of
$\delta^\ast_S$ projecting onto $h_0$. Only one
such plaquette, however, is $1$-connected with
the c-plaquettes. Since $\delta^\ast_S$ is $1$-connected, the
other plaquettes projecting onto $h_0$ must be  $1$-connected
with at least one other plaquette of
$\delta^\ast_S$. Each
of these further plaquettes projects into $\pi(C)$,
in contradiction of Lemma~\ref{properties}(iv).

We may now verify (\ref{Piproof1}) as follows.
Since $h_1$ is a w-plaquette, there exists $h_2\in A\cup B$, $h_2\ne h_1$,
such that $\pi(h_2)=h_0$. If there exists such $h_2$ belonging to
$A$, then (\ref{Piproof1}) holds. We assume the contrary, and let $h_2$ be such a plaquette
with $h_2\in B$. Since $h_1\in A$, for every $\eta\in U(h_0)\cap\pi(S)$,
$\eta\ne h_0$, there exists $\eta'\in A$ such that $\pi(\eta')\subseteq
[\eta]$ and $\eta'\1sim h_1$. [If this fails for some $\eta$, then,
as in the proof of Lemma \ref{properties}(ii),
in any configuration with interface $\de_S$,
there exists a path of open edges joining the vertex just above $h_1$
to the vertex just beneath $h_1$.  Since, by assumption,
all plaquettes of $A\cup B$ other than $h_1$, having projection $h_0$,
lie in $B$, this contradicts the fact that $\de_S$ is an interface.]
If any such $\eta'$ is vertical, then (\ref{Piproof1}) follows.
Assume that all such $\eta'$ are horizontal. Since $h_2\in B$,
there exists $h_3\in A$ such that $h_3\1sim h_2$, and
(\ref{Piproof1}) holds in this case also.

\noindent
(ii) The second part of the lemma follows from the observation that each of the plaquettes
in $A$ is $1$-connected to no more than four horizontal plaquettes of $B$.

\noindent
(iii)
Recall from  the remark after (\ref{hbnd}) that $A$ contains at least $|\pi(S)|$
horizontal plaquettes. Furthermore, $A$ must contain
at least $D(S)$ vertical plaquettes, and the claim follows.
\QED

Finally in this section, we derive an exponential
bound for the number of groups of walls satisfying certain constraints.


\begin{lem}
\label{numberofpif}
Let $h\in \delta_0$. There exists a constant $K$ such that\/\mbox{\rm:} 
the number of groups of
walls $G\in {\cal G}_{h}$ satisfying $\Pi(G)=k$ is no greater than $K^k$.

\end{lem}

\demo{Proof}
Let $G=\{S_1,S_2, \ldots, S_n\}\in {\cal G}_{h}$ where the $S_i=(A_i,B_i)$
are non-empty standard walls and $S_1\in {\cal S}_{h}$.
For $j\in\de_0$, define
$$
  R_j=\Big\{h'\in \delta_0: \|j,h'\|\leq \sqrt{\rho(j, \delta_G)} \Big\}\setminus \pi(G)
$$
and
$$
   \widetilde{G}=\Big(\bigcup_{i=1}^n [A_i\cup B_i]\Big) \cup \Big(\bigcup_{j\in \pi(G)}R_j\Big).
$$
There exist constants $C'$ and $C''$ such that, by  Lemma~\ref{Pi},
$$
 | \widetilde{G}|\leq |G|+C'\sum_{j\in \pi(G)}\rho(j, \delta_G)
\leq C'' |G|\leq 5\cdot 14 C'' \Pi(G),
$$
where $|G|=|\bigcup_i (A_i\cup B_i)|$.

It may be seen that $ \widetilde{G}$ is a
$0$-connected set of plaquettes containing $h$. Moreover, the $0$-connected sets
obtained by removing all the horizontal plaquettes $h'\in  \widetilde{G}$, for
which there exists no other plaquette $h''\in  \widetilde{G}$ with $\pi(h'')=\pi(h')$, are the standard walls
of $G$. Hence, the number of such groups of walls with $\Pi(G)=k$ is no greater than the number of
$0$-connected sets  of plaquettes containing no more than
$70 C'' k$ elements including $h$. It is proved in \cite{Do}, Lemma~2, that
there exists $\nu<\infty$ such that the number of $0$-connected sets of size $n$ containing $h$
is no larger than $\nu^n$. Corresponding to each such set there
are at most $2^n$ ways of partitioning the plaquettes between the $A_i$ and
the $B_i$. The claim of the lemma follows.
\QED


\section*{8. Exponential bounds for probabilities}
Let $\GG=\{G_h: h\in\de_0\}$ be a family of groups of walls.
If $\GG$ is admissible, there exists by 
Lemma \ref{construction} a unique
corresponding interface $\de_{\GG}$.
We may pick a random group $\zeta=\{\zeta_h: h\in\de_0\}$
of walls according to the
probability measure $\Prob_L$ induced by $\phi_L$ thus:
$$
\Prob_L(\zeta=\GG)
     =\left\{\begin{array}{ll}
  \ol\phi_{L}(\Delta=\de_{\GG})
         & \mbox{ if } \GG \mbox{ is admissible},\\
0 & \mbox{ otherwise.}
\end{array}\right.
$$



\begin{lem}
\label{probf}
Let $q\ge 1$, and let $p^*$ be as in Proposition \ref{propositionf}.
There exist constants $C_3$, $C_4$ such that
$$
\Prob_L\Big(\zeta_{h'}=G_{h'}\Bigmid\zeta_h=G_{h}
   \mbox{\rm\ for }h\in \delta_0,\ h\ne h'\Big)
       \leq 
C_3[C_4(1-p)]^{\Pi(G_{h'})},
$$
for $p>p^*$, and for all $h'\in \delta_0$, $G_{h'}\in {\cal G}_{h'}$, $L>0$, 
and for any
admissible family
$\{G_{h}: h\in \delta_0,\ h\ne h'\}$ of groups of walls.
\end{lem}

\demo{Proof}
The claim is trivial if $\GG=\{G_h: h\in\de_0\}$ is not admissible,
and therefore we may assume it admissible. Let $h'\in\de_0$,
and let $\GG'$ agree with
$\GG$ except at $h'$, where $G_{h'}$ is replaced by the 
empty group ${\cal E}_{h'}$. We write $\de=\de_\GG$ and
$\de'=\de_{\GG'}$.
Then
\be
\label{probfproof}
\Prob_L\Big(\zeta_{h'}=G_{h'}\Bigmid \zeta_h=G_h \mbox{ for }
    h\in\delta_0,\ h\ne h'\Big)
    \leq \frac{\ol\phi_L(\delta)}
                              {\ol\phi_L(\delta')}.	  
\ee
We will use (\ref{probdelta}) to bound the 
right-hand side of this expression.
In doing so, we shall require bounds for
$|\delta|-|\delta'|$, 
$|\odelta\setminus \delta|-|\odelta'\setminus\delta'|$, 
$K_{\delta}- K_{\delta'}$, and 
\be
\label{closeddif}
   \sum_{e\in E(\delta)\cap E_L} f_p(e,\delta, L)
                 -\sum_{e\in E(\delta')\cap E_L}f_p(e,\delta', L). 
\ee

It is easy to see from the definition of $\de$ that
$$
|\delta|=|\delta_0|+\sum_{h\in \delta_0}\Big[N(G_h)-|\pi(G_h)|\Big],
$$ 
and it follows that 
\be
\label{deldel}
|\delta|-|\delta'|= N(G_{h'})-|\pi(G_{h'})|= \Pi(G_{h'}).
\ee

A little thought leads to the inequality
\be
\label{ineq1}
|\odelta\setminus \delta|-|\odelta'\setminus\delta'|\geq 0,
\ee
and the reader may wish to omit the explanation which follows.
We claim that (\ref{ineq1}) follows from the inequality
\be
\label{ineq2}
|P(\odelta)|- |P(\odelta')|\ge 0,
\ee
where $P(\odelta)$ (respectively $P(\odelta')$) is
the set of plaquettes in
$\odelta\setminus\delta$ (respectively  $\odelta'\setminus\delta'$) 
which project into $[\pi(G_{h'})]$. In order to see
that (\ref{ineq2}) implies (\ref{ineq1}), we argue
as follows.
We may construct the extended interface $\odelta$ from $\odelta'$
in the following manner. First we
remove all the plaquettes from $\odelta'$ that project into $[\pi(G_{h'})]$, and
we fill the gaps by introducing the walls of $G_{h'}$ 
one by one along the lines of the 
proof of Lemma~\ref{construction}. Then we add the  
plaquettes of $\odelta\setminus \delta$ that project into 
$[\pi(G_{h'})]$. 
During this operation on interfaces,
we remove $P(\odelta')$ and add $P(\odelta)$; the claim follows. 

By Lemma~\ref{properties}(viii), there exists no vertical plaquette of 
$\odelta'\setminus \delta'$ that projects into $[\pi(G_{h'})]$ and is
in addition $1$-connected to some wall 
not belonging to $G_{h'}$. Moreover, since all the horizontal
plaquettes of $\odelta'$ belong to the semi-extended interface  
$\delta'^\ast$, those that
project onto $[\pi(G_{h'})]$ are c-plaquettes of $\delta'^\ast$; hence, 
such plaquettes lie in $\de'$. It follows that $P(\odelta')$ comprises the vertical
plaquettes that are $1$-connected with $\pi(G_{h'})$.

It is therefore sufficient to construct an injective map $T$ that maps each 
vertical plaquette $1$-connected with $\pi(G_{h'})$ to a different vertical plaquette in
$P(\odelta)$. 
We noted in the proof of Lemma \ref{Pi}(i)
that, for every $h_0\in \pi(G_{h}')$, there exists a horizontal plaquette 
$h_1\in \delta$ with $\pi(h_1)=h_0$.   
For every vertical
plaquette $h^\ve \1sim h_0$, there exists a translate
$h^\ve_1 \1sim h_1$. 
Suppose $h^\ve$ lies above $\delta_0$. If $h^\ve_1\in \odelta\setminus\de$, 
we set $T(h^\ve)=h^\ve_1$. If $h^\ve_1\in\de$, we consider the
(unique) vertical plaquette `above' it, which we denote $h^\ve_2$.  We repeat this procedure up to 
the first $n$ that  we meet a plaquette $h^\ve_n\in\odelta\setminus\de$,
and we set $T(h^\ve)=h^\ve_n$. When 
$h^\ve$ lies below $\delta_0$, we act similarly to find a plaquette $T(h^\ve)$
of $\odelta\setminus\de$ beneath $h^\ve$. The resulting $T$ is as required.

Turning to $K_\de-K_{\de'}$, we recall the notation after Proposition \ref{Q}.
Notice that exactly two of the components $(S_\de^i,U_\de^i)$
are infinite, and we suppose that these are assigned indices $1$ and $2$. 
For $i=3,4, \ldots, K_\delta$,  
let $ H(S_\delta^i)$ be the set of plaquettes that are the dual 
to an edge having
one vertex in $S_\delta^i$ and one vertex in $\partial S_\delta^i $. The finite component
$(S_\delta^i, U_\delta^i)$ is in a natural way surrounded by a particular wall,
namely that to which all the plaquettes of $H(S_\delta^i)$ belong. 
This follows from Lemma~\ref{properties}(v, viii) 
and the facts that 
$$
  P_i =\Big\{\pi\Big(h(\bond{x,x+(0,0,1)})\Big): x\in S_\delta^i \Big\}
$$   
is a $1$-connected subset of $\delta_0$, and that 
$[\pi(H(S_\delta^i))]=[P_i]$.   

Therefore,
\be
\label{kbnd0}
  K_\delta-K_{\delta'}= K_{\delta''}-2, 
\ee
where $\de''=\delta_{G_{h'}}$. It is elementary 
by Lemma \ref{Pi}(i) that
\be
\label{kbnd}
K_{\delta''}\leq 2N(G_{h'})\leq 28 \Pi(G_{h'}).
\ee

Finally, we estimate (\ref{closeddif}).
Let $H_1,H_2,\ldots, H_r$ be the maximal $0$-connected sets of plaquettes  in 
$\delta_0\setminus \pi(G_{h'})$, and let 
$\delta_i$ (respectively  $\delta'_i)$ be the set of plaquettes  of 
$\delta$ (respectively  $\delta'$) that project into $[H_i]$. Recalling the construction of 
an interface from its standard walls in the proof of Lemma \ref{construction}, 
there is a natural one--one correspondence between 
the plaquettes  of   $\delta_i$ and those of $\delta_i'$, and hence between the 
plaquettes in $U=\bigcup_{i=1}^r\delta_i$ and those in $U'=\bigcup_{i=1}^r \delta_i'$.
We denote by $T$ the  corresponding bijection that maps an edge $e$ with  
$h(e)\in \bigcup_{i=1}^r\delta_i$ to the edge $T(e)$ with corresponding dual plaquette in  
$\bigcup_{i=1}^r \delta_i'$. Note that $T(e)$ is a vertical translate of $e$.

If $e$ is such that $h(e)\in U$,
$$
  G(e, \delta, E_L ;T(e), \delta', E_L)
      \geq \|\pi'(h(e)),\pi(G_{h'})\| -1, 
$$
where $\pi'(h) $ is the earliest plaquette $h''$ of $\delta_0$ 
such that $\pi(h)\subseteq [h'']$, and
$$
\|h_1,H\| = \min\Big\{\|h_1,h_2\|: h_2\in H\Big\}.
$$
Let $p>p^\ast$.
Using the notation of Proposition \ref{propositionf} and Lemma \ref{limitf},
\begin{eqnarray}\label{bigbnd}
&&\kern-15pt\left|\sum_{e\in E(\delta)\cap E_L} f_p(e,\delta,L)-
     \sum_{e\in E(\delta')\cap E_L} f_p(e,\delta',L)\right|\nonumber\\
&&  \displaystyle     \hspace{.1cm} \leq \sum_{e\in E(U)\cap E_L}
                            |f_p(e,\delta, L)- f_p(T(e),\delta',L) |\nonumber\\
&&  \displaystyle     \hspace{3.5cm} +\sum_{e\in E(\delta\setminus U)\cap E_L}		    
                        f_p(e,\delta, L)
  +\sum_{e\in E(\delta'\setminus U')\cap E_L}   
                        f_p(e,\delta', L)\nonumber\\
&&  \displaystyle     \hspace{.1cm}\leq C_2 e^\gamma\sum_{e\in E(U)\cap E_L}
           \exp\Big(-\gamma\|\pi'(h(e)),\pi(G_{h'})\|\Big)
+C_1\Big[N(G_{h'})+|\pi(G_{h'})|\Big].\hspace{1cm}{}
\end{eqnarray}
By Lemma~\ref{Pi}, the second term of the last line 
is no greater than $C_5\Pi(G_{h'})$ for some constant
$C_5$.  Using the same lemma 
and the definition of a group of walls, we see that the first term is no larger than 
\begin{eqnarray}
\label{bnd8}
&&    C_2e^\gamma\sum_{h\in \delta_0\setminus \pi(G_{h'})}
           \rho(h, \delta)\exp\Big(-\gamma\| h ,\pi(G_{h'})\|\Big)\nonumber\\
&&      \hspace{2cm}\leq C_2e^\gamma\sum_{h\in \delta_0\setminus \pi(G_{h'})}
           \| h ,\pi(G_{h'})\|^2 \exp\Big(-\gamma\| h ,\pi(G_{h'})\|\Big)\nonumber\\	
&&     \hspace{2cm}\leq C_2e^\gamma\sum_{h''\in\pi(G_{h'})}\,
           \sum_{h\in \delta_0\setminus \pi(G_{h'})}
           \| h ,h''\|^2 \exp\Big(-\gamma\| h ,h''\|\Big)\nonumber\\
&& \hspace{2cm}\le  C_6|\pi(G_{h'})| \le 13C_6 \Pi(G_{h'}),
\end{eqnarray}
for some constant $C_6$.

The required conditional probability is, by (\ref{probdelta})
and (\ref{probfproof}),
$$
  p^{|\odelta\setminus\de|-|\odelta'\setminus\de'|} (1-p)^{|\de|-|\de'|}
q^{K_\de-K_{\de'}} \exp\left(\sum_{e\in E(\delta)\cap E_L} f_p(e,\delta,L)-
     \sum_{e\in E(\delta')\cap E_L} f_p(e,\delta',L)\right),
$$
which, by (\ref{deldel})--(\ref{bnd8}), is bounded above as in the statement of the lemma.
\QED

\section*{9. Main theorem}
Let $h\in\de_0$. For $\omega\in\Omega^\mu_{L}$, we write
$h \lra \infty$ if there exists a sequence $h=h_0, h_1,\dots,h_r$ 
of plaquettes in $\de_0$ such that:
$h_i\1sim h_{i+1}$ for $0\le i < r$; each $h_i$ is a c-plaquette of $\Delta(\omega)$;
$h_r=h(e)$ for some $e \notin E_L$.

\begin{theo}
\label{theorem1}
Let $q\geq 1$. For all $\epsilon>0$, there exists 
$\widehat p=\widehat p(\epsilon)<1$ such that, if $p>\widehat p$,  
\be
\label{theorem2res}
  \ol\phi_L(h \lra\infty)>1-\epsilon
\ee
for all $h\in \delta_0$ and all $L\ge 1$.
\end{theo}

Since, following Theorem \ref{theorem1}, $h$ is a c-plaquette with high probability,
it follows by Proposition \ref{Q} and the discussion immediately thereafter
that the vertex of $\ZZ$ immediately beneath (respectively above)
the centre of $h$ is joined to
$\partial^-\La_L$ (respectively $\partial^+\La_L$)
with high probability.
Thus Theorem \ref{theorem2} holds. Furthermore, since $h\lra\infty$ with
high probability, such connections
may be found within the plane of $\ZZ$ comprising vertices $x$ with $x_3=0$
(respectively $x_3=1$). 

 The existence of non-translation-invariant (conditioned) 
random-cluster measures
follows from Theorem \ref{theorem1},
as in the following sketch argument. 
For $e\in\Ed$, we write $e^\pm=e\pm(0,0,1)$.
Let $\omega\in\Omega$.
If $h=h(e)\in\de_0$ is a c-plaquette of $\Delta(\omega)$, then
$e$ is closed, and $h(e^\pm)\notin\ol{\Delta(\omega)}$. The 
configurations in the two regions
above and below $\Delta(\omega)$ are governed by wired random-cluster measures.
[We have used Lemma 8 here.] Hence, under (\ref{theorem2res}),
$$
\ol\phi_L(\omega(e)=1)\le\epsilon,\q \ol\phi_L(\omega(e^\pm)=1)\ge \frac{(1-\epsilon)p}
{p+(1-p)q},
$$
by Lemma \ref{comp}. Note that these inequalities concern the probabilities
of cylinder events.

Our second main result concerns the vertical 
displacement of the interface, and states 
roughly that there exists a geometric
bound on the tail of the displacement, uniformly in $L$.
Let $\de\in\DL$, $(x_1,x_2)\in\Z^2$, and write $x=(x_1,x_2,\tfrac12)$.
We define the {\it displacement\/} of $\de$ at $x$ by
$$
D(x,\de)=\sup\Bigl\{|d-\tfrac12|: (x_1,x_2,d)\in[\de]\Bigr\}.
$$

\begin{theo}
\label{theoremht}
Let $q\ge 1$. There exists $\widehat p <1$ and $\alpha(p)$ satisfying
$\alpha(p)>0$ when $p>\widehat p$ such that
$$
\ol\phi_L\Bigl(D(x,\Delta)\ge d\Bigr) \le e^{-d\alpha(p)} \quad
\mbox{for } d\ge 1,
$$
for all $(x_1,x_2)\in\Z^2$ and $L\ge 1$.
\end{theo}

\demo{Proof of Theorem \ref{theorem1}}
Let $h\in\de_0$.
We have not so far specified the ordering of 
plaquettes in $\delta_0$ used to 
identify the origin of a standard  wall or of a group of walls. 
We assume henceforth that 
this ordering is such that: for all $h_1,h_2\in\de_0$,
$h_1> h_2$ implies  $\|h,h_1\|\ge \|h,h_2\|$. 

For any standard wall $S$ there exists, by Lemma \ref{properties}(vi),
a unique maximal infinite 1-connected component $I(S)$
of $\de_0\setminus\pi(S)$.
Let $\omega\in\Omega^\mu_L$. The interface $\Delta(\omega)$
gives rise to a family of standard walls, and
$h\lra\infty$ if and only if, for each such wall $S$,
 $h$ belongs to $I(S)$. (This is a consequence
of a standard property of $\Z^2$;
see the appendix of \cite{Ke82}.)
Suppose on the contrary that $h\notin I(S_j)$
for some such standard wall $S_{j}$, for some $j\in\de_0$, 
belonging in turn to
some maximal admissible group $G_{h'}\in {\cal G}_{h'}$ of walls of $\Delta$,
for some $h' \in \delta_0$.
We have by  
Lemma \ref{Pi} and the above ordering on members of $\de_0$  that 
$$
  13 \Pi(G_{h'})\geq |\pi(G_{h'})|\geq |\pi(S_j)|\geq \|h,j\|+1\geq \|h,h'\|+1.
$$

Let $K$ be  as in Lemma~\ref{numberofpif}, and $p^*$, $C_4$  as in 
Lemma~\ref{probf}. We let $\widetilde p$ be sufficiently
large that $\widetilde p > p^*$ and 
$$
\lambda=\lambda(p)=-\tfrac1{13}\log[KC_4(1-p)]
$$
satisfies $\lambda(\widetilde p) > 0$.
By the latter lemma,
when $p>\widetilde p$,
\begin{eqnarray*}
1-\ol\phi_L(h \lra\infty)
    &\leq& \sum_{h'\in\delta_0}\Prob_L\Big(\Pi(\zeta_{h'})\geq \tfrac{1}{13}(\|h,h'\|+1)\Big)\\
    &\leq& \sum_{h'\in\delta_0}\, \sum_{n\geq (\|h,h'\|+1)/13}\,
         \sum_{{G\in {\cal G}_{h'}:}\atop{\Pi(G)=n}}\Prob_L(\zeta_{h'}=G)\\
    &\leq& \sum_{h'\in\delta_0}\, \sum_{n\geq (\|h,h'\|+1)/13} K^n C_3\big[ C_4(1-p)\big]^n\nonumber\\ 
    &\leq& C_3\sum_{h'\in \delta_0}\exp(-\lambda(\|h,h'\|+1))
    \leq C_7 e^{-\lambda},
\end{eqnarray*}
for appropriate constants $C_i$. The claim follows on choosing
$p$ sufficiently close to 1.
\QED

\demo{Proof of Theorem \ref{theoremht}}
This is related to the proof of Proposition 2.4 of \cite{BLP}.
If $D(x,\Delta)\ge d$, there exists $r$ satisfying $1\le r\le d$ such that
the following statement holds.
There exist distinct plaquettes $h_1,h_2,\dots,h_r\in\de_0$, and
maximal admissible groups $G_{h_i}$, $1\le i\le r$, of walls
of $\Delta$ such that: $x=(x_1,x_2,\tfrac12)$ lies in the interior of
one or more standard wall of each $G_{h_i}$, and $\sum_{i=1}^r \Pi(G_{h_i})
\ge d$ (recall Lemma \ref{Pi}(iii)).  
Let $m_i=\lfloor\tfrac1{13}(\|x,h_i\|+1)
\rfloor$
where $\|x,h\|=\|x-y\|$ and $y$ is the centre
of $h$.
By Lemma \ref{probf}, and as in the previous proof,
\begin{eqnarray*}
\ol\phi_L(D(x,\Delta)\ge d)
 &\le& \sum_{{h_1,h_2,\dots,h_r}\atop{1\le r\le d}}
 \Prob_L\left(\sum_i\Pi(\zeta_{h_i})\ge d,\ 
 \Pi(\zeta_{h_i})\ge  m_i\vee1\right)\\
&=& \sum_{{h_1,h_2,\dots,h_r}\atop{1\le r\le d}}\,\sum_{s=d}^\infty \,
\sum_{\stackrel{z_1,z_2,\dots,z_r:}{{z_1+z_2+\cdots+z_r=s}\atop{z_i\ge m_i\vee1}}}
\Prob_L\Big(\Pi(\zeta_{h_i})=z_i \mbox{ for } 1\le i \le r\Big)\\
&\le& \sum_{h_i}\sum_{s\ge d}  C_8[KC_4(1-p)]^s
\sum_{\stackrel{z_1,z_2,\dots,z_r:}{{z_1+z_2+\cdots+z_r=s}\atop{z_i\ge m_i\vee 1}}}
1,
\end{eqnarray*}
for some constant $C_8$.
The last summation is the number of ordered partitions of the integer
$s$ into $r$ parts, the $i$th of which is
at least $m_i\vee 1$. By adapting the
classical solution to this enumeration
valid for the case $m_i\equiv 1$ (see, for example, \cite{Hall}),
we see that
$$
\sum_{\stackrel{z_1,z_2,\dots,z_r:}{{z_1+z_2+\cdots+z_r=s}\atop{z_i\ge m_i\vee 1}}}
1 \le {{s-1-\sum_i m_i\vee 1}\choose{r-1}}\le
2^{s-1-\sum_i m_i\vee 1}\le 2^{s-1-\sum_i m_i},
$$
whence, for some $C_9$,
$$
\ol\phi_L(D(x,\Delta)\ge d)
 \le C_9\sum_{s\ge d} [2KC_4 (1-p)]^s
\left[\sum_{h\in\de_0} 2^{-\lfloor\|x,h\|/13\rfloor}\right]^d,
$$
which decays exponentially as $d\to\infty$
when $2KC_4(1-p)$ is sufficiently small.
\QED

\section*{Acknowledgement}
We thank Roman Koteck\'y for his comments on an aspect of this work.

\end{document}